\documentclass[12pt,reqno,dvips]{amsart}
\usepackage{amsmath,amsthm,amsfonts,amssymb,amscd,amstext}
\usepackage[ansinew]{inputenc}
\usepackage{graphicx}
\usepackage{psfrag,euscript}
\usepackage{a4wide}
\usepackage{pstricks}
\usepackage{mathpazo} 

\numberwithin{equation}{section}

\usepackage{color}
\usepackage[colorlinks=true,linkcolor=red,urlcolor=black]{hyperref} 

\newcommand{\qand}{\quad\text{and}\quad}

\theoremstyle{plain}
\newtheorem{maintheorem}{Theorem}
\newtheorem{maincorollary}[maintheorem]{Corollary}

\newtheorem{theorem}{Theorem}[section]
\newtheorem{proposition}[theorem]{Proposition}

\newtheorem{lemma}[theorem]{Lemma}

\newtheorem{problem}{Problem}
\theoremstyle{definition}
\newtheorem{remark}[theorem]{Remark}

\newcommand{\RR}{{\mathbb R}}

\newcommand{\ZZ}{{\mathbb Z}}

\newcommand{\vfi}{\varphi}

\newcommand{\de}{\delta}

\renewcommand{\epsilon}{\varepsilon}

\newcommand{\supp}{\operatorname{supp}}

\newcommand{\cF}{\EuScript{F}}

\newcommand{\F}{\EuScript{F}}

\newcommand{\V}{\EuScript{V}}
\newcommand{\U}{\EuScript{U}}

\newcommand{\Diff}{{\rm Diff}^r(M)}

\newcommand{\ecu}{E^{cu}}

\newcommand{\m}{{\rm Leb}\:}

\title[Abundance of non-zero central Lyapunov exponents] {On
  the abundance of non-zero central Lyapunov exponents,
  physical measures and stable ergodicity for partially
  hyperbolic dynamics}

\begin{thanks} {V.A. was partially supported by CNPq, FAPERJ
    and PRONEX (Brazil). C.V. was partially supported by Proyecto Fondecyt
    1100547 and Research Network on Low Dimensional
    Dynamics, PBCT ACT-17, CONICYT, Chile.  }
\end{thanks}

\author{V\'{\i}tor Ara\'ujo}

\address{V\'\i tor Ara\'ujo, Instituto de Matem\'a\-tica,
  Universidade Federal do Rio de Janeiro, C. P. 68.530,
  21.945-970 Rio de Janeiro, RJ-Brazil.} 
\email{vitor.araujo@im.ufrj.br \textrm{or} vdaraujo@impa.br}

\author{Carlos H. V\'asquez}

\address{Carlos H. V\'asquez, Instituto de Matem\'atica,
Pontificia Universidad Cat\'olica de Valpara\'{\i}so,
Blanco Viel 596, Cerro Bar\'on, 
Valpara\'{\i}so-Chile.}
\email{carlos.vasquez@ucv.cl}

\date{\today}

\begin{document}


\renewcommand{\subjclassname}{\textup{2000} Mathematics Subject Classification}


\begin{abstract}
  We show that the time-1 map of an Anosov flow, whose
  strong-unstable foliation is $C^2$ smooth and minimal, is
  $C^2$ close to a diffeomorphism having positive central
  Lyapunov exponent Lebesgue almost everywhere and a unique
  physical measure with full basin, which is $C^r$ stably
  ergodic.  Our method is perturbative and does not rely on
  preservation of a smooth measure.

\bigskip
\noindent
\textsc{Résumé}:
  Nous montrons que le temps-$1$ d'un flux d'Anosov, dont
  le foliation forte-instable est $C^2$ lisse et minimal,
  est $C^2$ proche d'un difféomorphisme ayant exposant de
  Lyapunov central positif Lebesgue presque partout et une
  unique mesure physique avec bassin plein, ce qui est
  $C^r$-stablement ergodique. Notre méthode est
  perturbatif, et ne repose pas sur la préservation d'une
  mesure de volume.
\end{abstract}

\maketitle

\tableofcontents

\section{Introduction}
\label{sec:introd}

Stable ergodicity is a desirable property for dynamical
systems since it is arguably the most basic global
statistical feature. This is inspired by the fundamental
Boltzman Ergodic Hypothesis from Statistical Mechanics,
which is the main motivation behind the celebrated Birkhoff
Ergodic Theorem, ensuring the equality between temporal and
spatial averages with respect to a (ergodic) probability
measure $\mu$ invariant under a measurable transformation
$f:M\to M$ of a compact manifold $M$, i.e. for every
integrable function $\vfi:M\to\RR$ we have
\begin{align}
  \label{e-birkhoff}
  \lim_{n\to+\infty}\frac1n\sum_{j=0}^{n-1}
  \vfi\big( f^j(x) \big) = \int \vfi \, d\mu
\end{align}
for $\mu$ almost every point $x\in M$.

Systems nearby a stably ergodic system remain ergodic. In
general we need some ``natural measure'' to explore stable
ergodicity. D. V. Anosov in \cite{An67} was the first to
establish the existence of open sets of ergodic systems on a
wide class of manifolds: the geodesic flow on the unit
tangent bundle of compact Riemannian manifolds with constant
negative sectional curvature. Systems sharing the same
features are known today as ``Anosov systems'': they are
globally hyperbolic and structurally stable (all dynamical
properties of their perturbations are the same on all
scales). Previous results on hyperbolicity and ergodicity
with respect to the natural Liouville volume measure on the
unit tangent bundle of compact surfaces with constant negative
Gaussian curvature were obtained earlier by
Hedlund~\cite{Hedl34} and Hopf~\cite{hopf39,Hopf71}.

Pugh and Shub began a program devoted to stable ergodicity
\cite{PuSh2004} for partially hyperbolic volume preserving
systems. They ask how frequently is partial hyperbolicity
the main reason for a dynamical system to be stably ergodic
and conjectured that, among the volume preserving partially
hyperbolic dynamical systems, the stably ergodic ones form
an open and dense set. Many advances have been obtained in
this direction recently; see e.g. Rodriguez Hertz et
al~\cite{HHertzUres07}. The natural measure in this setting
is the Lebesgue volume measure which disintegrates as a
density along the unstable foliations associated to partial
hyperbolic dynamics: such measures are known as
\emph{u-Gibbs states}.

In the dissipative setting the \emph{SRB/physical measures}
are natural candidates to play the role of Lebesgue measure,
since for these measures the time averages coincide with the
space averages for orbits starting in a positive volume
subset of the ambient space: points satisfying
\eqref{e-birkhoff} for all continuous functions
$\vfi:M\to\RR$ (the \emph{ergodic basin} of the measure)
form a subset of positive volume. Moreover, such physical
measures for partially hyperbolic systems disintegrate as
densities along the unstable directions also, that is, they
are $u$-Gibbs states. So the question raised by Pugh
and Shub makes sense also in a non-conservative setting.

In \cite{BuDoPe2002,BuDoPePo2008} Burns, Dolgopyat, Pesin
and Pollicott studied stable ergodicity for partially
hyperbolic diffeomorphims, not necessarily conservative,
whose Lyapunov exponents along the center direction are all
negative with respect to some $u$-Gibbs state which, in this
setting, becomes a physical measure. In \cite{vasquez2009,
  AndVas2011} the case with positive Lyapunov exponents
along the center direction was considered.  The main reason
to consider non-zero central Lyapunov exponents in the works
cited before is that in many cases hyperbolic $u$-Gibbs
states are physical measures \cite{Do2000, BoV00,
  ABV00,vasquez2009}. Moreover Bochi, Fayad and Pujals in
\cite{BoFaPuj} have noted that stably ergodic
\emph{conservative} diffeomorphisms must be $C^1$ close to a
conservative diffeomorphism with stably non-zero central
Lyapunov exponents.

We recall the following problem posed by Bonatti, D\'iaz and
Pujals in \cite{BDP}.

\begin{problem}\label{prob:1}
  Let $f$ be a $C^1$ robustly transitive diffeomorphism of
  class $C^2$ on a compact manifold $M$. Does there exist
  $g$ close to $f$ having finitely many physical measures
  such that the union of their basins has total Lebesgue
  measure in $M$?
\end{problem}

This problem is connected with partial hyperbolicity since
in \cite{BDP} it was proved that robustly transitive
diffeomorphisms in compact manifolds of dimension three are
partially hyperbolic.

Hence the abundance of non-zero Lyapunov exponents along the
central direction of a strongly partially hyperbolic
diffeomorphism is closely related with the issue of the
existence (abundance) of physical measures and its
ergodicity. This work is devoted to shed some light on this
subject.

\emph{The aim of our work is to show that, under certain
  conditions explained below, it is possible to remove zero
  central Lyapunov exponents by perturbation inside a class
  of partially hyperbolic diffeomorphisms which are not
  conservative.} This is an extension of results
from Burns, Pugh and Wilkinson\cite{Wi98,BPW00} about stable
ergodicity of the time-one map of a geodesic flow, now
encompassing also non-conservative perturbations.

There are several results about removing zero Lyapunov
exponents either in the conservative setting, by Shub and
Wilkinson \cite{ShWi00}, Ruelle \cite{Ruelle2003b} and
Baraviera and Bonatti \cite{BB03} for diffeomorphims and
\cite{BesRoc07} for flows; or for partially hyperbolic
diffeomorphisms with center foliation formed by compact
center leaves, but relying on rigidity arguments, by Viana
and Yang \cite{viana-yang} and F. Rodriguez-Hertz,
M.A. Rodriguez-Hertz, Tahzibi and Ures \cite{tahz010}.  On
the opposite direction, that of forcing zero Lyapunov
exponents in the absence of weak forms of hyperbolicity
(like dominated splitting on the tangent bundle) among
conservative diffeomorphisms, there are the works of
Bochi~\cite{Bochi02} together with Viana~\cite{BoV02}, and
the corresponding versions for flows by Bessa~\cite{bessa}
together with one of the authors \cite{arbes2008}, always
restricted to $C^1$ perturbation techniques.  However, to
the authors best knowledge, there are no other general
perturbation results on central Lyapunov exponents for
\emph{not conservative} systems.

We consider the case of the time-$1$ map of a $C^2$ Anosov
flow which is not a suspension flow of an Anosov
diffeomorphism and whose strong-unstable foliation is
smooth; see the next section for details.

Our method relies on the following motivation: if we assume
that a function $\vfi_f:M\to\RR$ and a $f$-invariant
probability measure $\mu_f$ are given, the function
depending continuously on the local dynamics of $f$ and the
measure depending continuously on $f$ in the space of
diffeomorphisms, then we should be able to change the value
of the integral $\int\vfi_f\,d\mu_f$ by an arbitrarily small
perturbation of the map $f$.

Our results are deduced applying the above idea to the
integral of the logarithm of the central Jacobian with
respect to a $u$-Gibbs state, for a strongly partially
hyperbolic system whose strong-unstable foliation is smooth
(the foliation is of class $C^2$) and minimal (every leaf is
dense in the ambient space).  By Oseledet's Multiplicative
Ergodic Theorem (see e.g.  Barreira
and Pesin~\cite{BarPes2002}) this integral gives the central
Lyapunov exponent. Therefore it is not possible to have a
constant zero central Lyapunov exponent for all $u$-Gibbs states
in a $C^2$ neighborhood of a partially hyperbolic map with a
smooth minimal strong-unstable foliation.

The continuity of $u$-Gibbs states under perturbations has
been studied by one the authors in \cite{vasquez2006}; see
also \cite{BuDoPePo2008}. To be able to control certain
features of the $u$-Gibbs states of the perturbed map, we
need to ensure that after the perturbation the new map is
$C^2$ close to the original one.  For this we were led to
assume that the strong-unstable foliation is $C^2$
smooth. In this setting, we can construct our perturbation
$C^2$ close to the original map, enabling us to use the
results of smooth ergodic theory already known for $u$-Gibbs
states.

The assumption of minimality is rather natural in this
setting since Parry, in \cite{parry70}, showed that a linear
torus automorphism is ergodic (with respect to the
Haar/volume measure) if, and only if, the corresponding
strong-stable foliation is minimal. Moreover by the results
of Plante~\cite{Plante72} for a transitive Anosov flow on a
compact manifold, either the strong-stable and
strong-unstable foliations are minimal, or the flow is the
suspension of an Anosov diffeomorphism of a compact
submanifold with codimension one.  We note that the
existence of partially hyperbolic diffeomorphisms having
robustly minimal strong-unstable foliations was obtained by
Bonatti, D\'iaz and Ures in \cite{BoDiUr}.

The minimality of the unstable foliation ensures, in our
setting, that \emph{the future orbit of Lebesgue almost
  every point has positive frequency of visits to any open
  subset.}  In the conservative setting, we have that
Lebesgue almost every point has well defined time averages
\emph{for the future and for the past}, i.e. under iterates
of the map and of the inverse map. This is not necessarily
true in general and prevents us from using arguments similar
to the ones of e.g. Burns, Wilkinson \cite{BuWi08} and
F. Rodriguez Hertz, M. A. Rodriguez Hertz, A. Tahzibi,
R. Ures \cite{HHTah07} in the conservative setting. 


The closely related results of Dolgopyat~\cite{Do2004} where
obtained using completely different techniques of a more
analytic nature, and are put in the setting of perturbation
along generic one-parameter families of maps through the
original map $f$. In addition the results are stated and
proved in the setting where the stable and unstable
subbundles are one-dimensional.

So, in this work we establish a perturbative method to
remove zero Lyapunov exponents (Theorem \ref{mthm:removing}
and its proof) in a non-conservative setting;
we deduce existence and uniqueness of physical measures
(Corollary~\ref{mthm:genericphysical}) and also
stable ergodicity for the perturbed systems (Corollary
\ref{mcor:openGibbs}).



\section{Statement of the results}
\label{sec:statem-results}

Let $M$ be a closed Riemannian manifold. We denote by
$\|\cdot\|$ the norm obtained from the Rimannian structure
and by ${\rm Leb\:}$ the Lebesgue measure on $M$.

If $V$, $W$ are normed linear spaces, we
define $$\|A\|=\sup\big\{\|Av\|_W/\|v\|_V, v\in
V\setminus\{0\}\big\},$$ and $${\mathrm
  m}(A)=\inf\big\{\|Av\|_W/\|v\|_V, v\in
V\setminus\{0\}\big\}$$ for a linear map $A:V\to W$.

A diffeomorphism $f\colon M\rightarrow M$ is {\em strongly
  partially hyperbolic} if there exists a continuous
$Df$-invariant splitting of $TM$, $$TM=E^s\oplus E^c\oplus
E^u,$$ and there exist constants $C\geq 0$ and

$$0<\lambda_1\leq \mu_1<\lambda_2\leq\mu_2 <\lambda_3\leq\mu_3$$ 
with $\mu_1<1<\lambda_3$ such that for all $x\in M$ and
every $n\geq 1$ we have:

\begin{equation}\label{PH1}
C^{-1}\lambda_1^n\leq {\rm m}\:(Df^n(x)|E^s(x))\leq\|Df^n(x)|E^s(x)\|\leq C\mu_1^n,
\end{equation}

\begin{equation}\label{PH2}
C^{-1}\lambda_2^n\leq {\rm m}\:(Df^n(x)|E^c(x))\leq\|Df^n(x)|E^c(x)\|\leq C\mu_2^n, 
\end{equation}

\begin{equation}\label{PH3}
C^{-1}\lambda_3^n\leq {\rm m}\:(Df^n(x)|E^u(x))\leq\|Df^n(x)|E^u(x)\|\leq C\mu_3^n. 
\end{equation}

%
%

The expression (\ref{PH1}) means that $E^s$ is uniformly
contracting, while (\ref{PH3}) means that $E^u$ is uniformly
expanding. The expression (\ref{PH2})  implies that $E^u$
dominates $E^c$ and that $E^c$ dominates $E^s$. We assume
that the subbundles are non-trivial.

It is well known that partially hyperbolic diffeomorphisms
have an unstable foliation $\cF^{u}= \{\cF^{u}(x)\::\:x\in
M\}$, whose leaves are the (strong) unstable manifolds
$W^u(x)$, $x\in M$; and 
they also have a stable foliation $\cF^s=\{
\cF^s(x)\::\:x\in M\}$, whose leaves are the (strong) stable
manifolds $W^s(x)$, $x\in M$. For $C^r$ diffeomorphisms,
$r>1$, these foliations are absolutely continuous; see
Hirsch, Pugh and Shub \cite{HPS77}.

An $f$-invariant probability measure $\mu$ is a {\em
  $u$-Gibbs state} if the conditional measures of $\mu$ with
respect to the partition into local strong-unstable
manifolds are absolutely continuous with respect to Lebesgue
measure along the corresponding local strong-unstable
manifolds.  If $f$ is a $C^2$-partially hyperbolic
diffeomorphism, there always exists a $u$-Gibbs state; see
Proposition~\ref{pr:uM1} or see e.g.~Pesin and Sinai~\cite{PS82} 
for more details.

We recall that,
if $\mu$ is a $f$-invariant measure, the (ergodic)
{\em basin} of $\mu$ is the set $B(\mu)$ of all points $x\in
M$ such that \eqref{e-birkhoff} is satisfied for all  continuous functions
$\vfi:M\to\RR$.
It is well known that the set $B(\mu)$ has full measure with
respect to any ergodic $f$-invariant probability measure
$\mu$. The $f$-invariant measure $\mu$ is a {\em physical}
or {\em SRB (Sinai-Ruelle-Bowen) measure}, if its basin
$B(\mu)$ has positive Lebesgue measure (volume) on $M$.


\subsection{Standing assumptions}
\label{sec:standing-assumpt}

We assume throughout that \emph{the unstable foliation
  $\cF^u$ is minimal}, that is, every leaf $\xi\in\cF^u$ is
dense in $M$. This property is satisfied by a $C^1$ open set
of partially hyperbolic diffeomorphisms by the results of
Bonatti, D\'\i az and Ures in \cite{BoDiUr}, among robustly
trasitive strongly partially hyperbolic diffeomorphism in
dimension three, that is, each subbundle is
one-dimensional. It is  satisfied by the time one map of
any transitive Anosov flow which is not the suspension of an
Anosov diffeomorphism of a codimension one submanifold, by
the results of Plante~\cite{Plante72}, e.g.,
the geodesic flow on surfaces of constant negative curvature
and by many contact Anosov flows.

This is important to ensure that certain properties of the
map obtained after local pertubations are spread to the
entire ambient space. 

We also assume that the subbundle $E^u$ induces a smooth
foliation $\cF^u$ of class $C^2$.  We remark that for the
general partially hyperbolic diffeomorphism the stable and
unstable laminations $\cF^s, \cF^u$, although having leaves
as smooth as $f$, are not foliations in the usual sense of
Differential Topology: their leaves do not ``stack on top of
each other'' in a smooth way.


\subsection{Removing zero central Lyapunov exponent}
\label{sec:pert-zero-central-0}

This is our main result. Abundance of non-zero central
Lyapunov exponents means the existence open sets of
diffeomorphisms where each one exhibits non-zero central
Lyapunov exponents for Lebesgue almost every point on the
manifold.
 
A partially hyperbolic
diffeomorphism such that every $u$-Gibbs state has positive
central Lyapunov exponents is called {\em mostly expanding}
and their properties are studied in
\cite{AndVas2011,vasquez2009}. Mostly expanding is the dual notion of {\em mostly contracting}
introduced by \cite{BoV00,Do2000} and  studied in
\cite{And2010}.

For $r\geq 2$, a $C^r$ mostly expanding diffeomorphism has
non-zero Lyapunov exponents Lebesgue almost everywhere in
the ambient manifold. This is a $C^r$-open property  and it
also implies the existence of physical measures  (see Proposition~\ref{mteo:MACA}).

\begin{maintheorem}\label{mthm:removing}
  Let $f$ be the time-$1$ map of an Anosov flow whose
  strong-unstable foliation is $C^2$ smooth and minimal. For
  $r\ge2$, $f$ is
  $C^2$-close to a $C^r$-open set of mostly expanding
  diffeomorphisms.
\end{maintheorem}

We give a brief sketch of our arguments.  Let $f$ be as in
the statement of Theorem~\ref{mthm:removing}.  In what
follows, for each $u$-Gibbs state $\mu$ of $f$, we
denote $$\lambda_\mu^c(f):=\int\log\|Df\mid E^c\| \,d\mu.$$
We note that, since the central direction is assumed to be
one-dimensional, if $\mu$ is ergodic, then this number
equals the central Lyapunov exponent. Moreover, since $f$ is
the time-$1$ map of an Anosov flow, we always have
$\lambda_\mu^c(f)=0$ for each $u$-Gibbs state $\mu$.


Using the perturbative methods explained in
Sections~\ref{sec:overvi-arguments},
\ref{sec:pert-centr-lyap}, \ref{sec:perturb-central-subb}
and~\ref{sec:compar-action-deriva}, we conclude that there
exist a partially hyperbolic diffeomorphism $g$ $C^2$-close
to $f$ such that $g$ is mostly expanding: for each $u$-Gibbs
state $\mu_g$ of $g$ we have $\lambda_{\mu_g}^c(g)>0$.  The
minimality of $\cF^u$ now ensures that there is a unique $\mu_g$ which is a
$cu$-Gibbs state and the unique physical measure; see
Lemma~\ref{le:unque-cuGibbs} for more details. The rest of
the conclusion follows from the $C^r$ openness of the mostly
expanding property (see Proposition~\ref{mteo:MACA} and Proposition~\ref{pr:stbl-erg-most-contr}).

More precisely, we show that, for $f\in \Diff$, $r\geq2$ in
the setting of Theorem~\ref{mthm:removing}, then \emph{
  arbitrarily $C^2$-close of $f$, there exists a
  $C^r$-neighborhood $\V$ of mostly expanding
  diffeomorphisms. As consequence, for each $g\in \V$ there
  exists a physical measure given by an ergodic $u$-Gibbs
  state $\mu$ for $g$ with positive central exponent, whose
  basin has full measure.}

\subsection{Abundance of physical measures with non-zero
  central exponent}
\label{sec:non-zero-central-1}

As a consequence of Theorem~\ref{mthm:removing} we provide a partial answer to
Problem~\ref{prob:1}.
As usual in smooth ergodic theory, we say that an invariant
probability measure $\mu$ is \emph{hyperbolic} if the
Lyapunov exponents of $\mu$-almost every point are never
zero.

\begin{maincorollary}\label{mthm:genericphysical}
  Each time-$1$ map of an Anosov flow, whose strong-unstable
  foliation is $C^2$ smooth and minimal, is $C^2$-close to a
  $C^r$-open set of partially hyperbolic diffeomorphims
  admitting a unique physical and hyperbolic measure with
  full basin (with $r\ge2$).
\end{maincorollary}

Hence Corollary~\ref{mthm:genericphysical} ensures that
Problem~\ref{prob:1} has an affirmative answer for time-$1$
maps $f$ of an Anosov flow when the strong-unstable foliation of
$f$ is $C^2$ smooth and minimal.

\subsection{Abundance of stable ergodicity}
\label{sec:statist-stochast-sta}

Now we rewrite our results from the point-of-view of stable
ergodicity.  A diffeomorphism $f$ is {\em $C^r$-stably
  ergodic} if there exists a $C^r$-neighborhood $\V$ of $f$,
where for each $g\in \V$ there exists a unique physical
measure $\mu_g$ whose ergodic basin has full Lebesgue
measure in the ambient space.

From the results of Section~\ref{sec:minimal-unstable-fol-1}
on uniqueness of physical measures for partially hyperbolic
diffeomorphims with minimal strong-unstable foliation, we
obtain the following.

\begin{maincorollary}
  \label{mcor:openGibbs}
  Let $f$ be the time-$1$ map of an Anosov flow whose
  strong-unstable foliation is $C^2$ smooth and minimal.Then
  $f$ is $C^2$-close to a $C^r$-stably ergodic
  diffeomorphism.
\end{maincorollary}

From the results \cite{BuDoPePo2008,vasquez2006} about the
continuous variation of $u$-Gibbs states and their
(non-zero) Lyapunov exponents with the diffeomorphism in the
$C^2$ topology, the stably ergodic diffeomorphisms we obtain
are necessarily statistically stable. This means that the
physical measures depend continuously on the diffeomorphism.


We now remark that Bochi, Fayad and Pujals in \cite{BoFaPuj}
have noted that stably ergodic \emph{conservative}
diffeomorphisms must be $C^1$ close to a conservative
diffeomorphism for which Lebesgue measure is ergodic and
hyperbolic, that is, stable ergodicity in the conservative
setting implies stably non-zero central Lyapunov exponents.

In our setting a straighforward consequence of our main
result reads as follows.

\begin{maincorollary}
  \label{mcor:staberg-zeroexp}
  The set of $C^2$ stably ergodic strongly partially
  hyperbolic diffeomorphisms with one-dimensional central
  subbundle and stably zero central Lyapunov exponents,
  cannot contain a diffeomorphism which is the time-$1$ map
  of an Anosov flow whose strong-unstable foliation is $C^2$
  smooth and minimal.
\end{maincorollary}


\subsection{Related open questions}
\label{sec:related-open-questi}

The results stated above suggest naturally the following
questions/problems.

\begin{problem}
  \label{prob:minimalityfoliation}
  We have seen that minimal unstable foliations are helpful
  to get a physical measure for partially hyperbolic
  systems. How general or abundant are the partially
  hyperbolic systems with minimal unstable foliations?
\end{problem}

\begin{problem}
  \label{prob:mixcentralbehavior}
  Given that we have a positive or negative central Lyapunov
  exponent for some $u$-Gibbs measure, for partially
  hyperbolic diffeomorphism with minimal unstable
  foliations, can we ensure that we have \emph{mixed central
    behavior} for higher dimensional central subbundles?
  That is, can we obtain generically or densely that there
  are hyperbolic $u$-Gibbs states with higher dimensional
  central subbundle having only non-zero Lyapunov exponents
  along this central direction?
\end{problem}  

\begin{problem}\label{prob:hyperbolicphysical}
  A natural problem in our setting is to understand if mixed
  central behavior together with hyperbolicity (absence of
  zero Lyapunov exponents along the central direction) for
  some $u$-Gibbs state is sufficient in general to obtain a
  physical measure. If not, what extra conditions are needed
  to obtain a physical measure in this setting?
\end{problem}

The last Corollary~\ref{mcor:staberg-zeroexp} naturally
suggests the following

\begin{problem}
  \label{prob:zerocentral}
  Is it possible to have stable ergodicity with zero
  Lyapunov exponents robustly along the center direction?
  That is, can we have stable ergodicity and zero central
  Lyapunov exponents almost everywhere simultaneously?
\end{problem}

\begin{problem}
  \label{prob:susceptibility}
  Since for our stably ergodic maps the physical measures vary continuously with the diffeomorphism,
is it true that the physical measures also depend \emph{smoothly} on the diffeomorphism? That is, can we
obtain a susceptibility function for strongly partially hyperbolic diffeomorphisms in our setting along
the lines of the work of Ruelle \cite{Ruelle2008,ruelle09}?
\end{problem}

The choice of the adequate $C^r$ topologies is part of the problems above.


\subsection{Overview of the arguments.}
\label{sec:overvi-arguments}

Here we present an overview of the arguments to be detailed
in what follows.

The statements of known results about $u$-Gibbs and
$cu$-Gibbs states versus physical measures, together with
minimality of strong-unstable foliation versus uniqueness of
$u$-Gibbs states are collected in
Section~\ref{sec:prelim-results} for convenience. We refer
to them along the rest of this text when needed.

  Let $f$ be the time-$1$ map of an Anosov flow whose
  strong-unstable foliation is $C^2$ smooth and minimal.
\emph{We claim that we can perturb $f$ to a $C^2$ close
  mostly expanding map $g$ whose corresponding
  strong-unstable foliation $\cF_g^u$ equals that of $f$:
  $\cF^u_g=\cF^u_f$.}

If we assume this claim, then since $\cF^u_g$ is
minimal we get from Lemma~\ref{le:unque-cuGibbs} that there
exists a unique $cu$-Gibbs state $\mu$ for $g$. Now we are
in the setting of
Proposition~\ref{pr:stbl-erg-most-contr}, and so we conclude
that there exists a $C^r$ neighborhood $\V$ of $g$ where all
maps are mostly expanding and have a unique physical measure
given by an ergodic  $cu$-Gibbs state.

This completes the proof of Theorem~\ref{mthm:removing} and
Corollary~\ref{mthm:genericphysical}, and also shows that
$g$ is $C^r$ stably ergodic as in the statement of
Corollary~\ref{mcor:openGibbs}, after we prove the claim above.

\subsubsection{Strategy}
\label{sec:strategy}

To prove the claim, we perform a local perturbation of the
map $f$ to a map $g$ by defining $g=f\circ H$, where $H$ is
a $C^2$ diffeomorphism of $M$ such that, for some given
non-periodic point $q_0\in M$ for $f$ and sufficiently small
$\epsilon,t>0$, we have in chosen local coordinates
\begin{itemize}
\item $H(B(q_0,2\epsilon))=B(q_0,2\epsilon)$;
\item $H$ is the identity map $Id$ on $M\setminus
  B(q_0,2\epsilon)$;
\item $\|H-Id\|_{C^1}=t\epsilon$ and
  $\|H-Id\|_{C^2}\xrightarrow[\epsilon\to0^+]{}0$;
\item $DH(q)\cdot E^c_f(q)$ is the graph of a
  \emph{non-zero} injective linear map $L_{H(q)}:E^c_f(H(q))\to
  E^u_f(H(q))$, for all $q\in
  B(q_0,2\epsilon)$;
\end{itemize}
where we write $E^*_f$ for the $Df$-invariant subbundles of
$f$, $*=s,c,u$.

We present the construction of this map $H$ in
Section~\ref{sec:pert-centr-lyap}, where the assumption of
$C^2$ smoothness on $\cF^u$ enables us to \emph{define $H$
  using coordinates along the leaves of the strong-unstable
  foliation and, most useful, to keep the strong-unstable
  and center-unstable foliations unchanged, so that $\cF^u$
  remains a minimal foliation for the perturbed map}.  This
perturbation $g=f\circ H$ of $f$, because it is $C^2$ close
to $f$, is also a strongly partially hyperbolic
diffeomorphism with $Dg$-invariant subbundles $E^*_g$,
$*=s,c,u$, where $E^c_g$ is one-dimensional.

We claim that the center $Dg$-invariant subbundle is
``tilted'' towards the original $E^u_f$ subbundle in such a
way that there exists a non-negative measurable function
$\xi\::\: M\to\RR$ such that for every $q\in M$,
\begin{equation}\label{eq:qoutient-gt-1}\|Dg(q)|E^c_g(q)\|\geq 1+\xi(q).\end{equation}
Moreover, the set of points $q\in M$ such that $\xi(q)>0$ has positive $\mu_g$ measure for every $u$-Gibbs state $\mu_g$ of $g$ (see Lemma~\ref{le:ML}). Then, we conclude that

$$
  \lambda_{\mu_g}^c(g)=
  \int\log\|Dg\mid  E^c_g\|\,d\mu_g \ge
  \int\log(1+\xi(q)) d\mu_g(q) >0.$$

Here $\xi(q)=\xi_{t,\epsilon}(q)\geq 0$  depends
on the domination of the action of $Df$ on $E^u_f$ over
$E^c_f$  (given by~\eqref{PH1}, \eqref{PH2} and
\eqref{PH3}) and on the $C^1$ distance between $f$ and $g$.
Both claims above are proved in
Sections~\ref{sec:perturb-central-subb}
and~\ref{sec:compar-action-deriva} by:
\begin{enumerate}
\item[(i)]  showing that the perturbed central subbundle
  $E^c_g$ has a non-zero component along the old unstable
  subbundle, in Section~\ref{sec:perturb-central-subb};
\item[(ii)] taking a Riemannian adapted norm $\|\cdot\|$ for
  the strongly partially hyperbolic diffeomorphims $f$,
  given by \cite{Goum07}, to estimate the expansion along
  $E^c_g$ in a transparent way, in
  Section~\ref{sec:compar-action-deriva}.
\end{enumerate}
The proof is mostly a linear algebra argument taking
advantage of the robustness of the domination of the
splitting.
These are the main arguments in the proof of
Theorem~\ref{mthm:removing} and
Corollaries~\ref{mthm:genericphysical} through
\ref{mcor:openGibbs}.


\subsection{Examples of application}
\label{sec:exampl-applic}

\subsubsection{The time-one map of the geodesic flow on
  surfaces of constant negative curvature}
\label{sec:time-one-map}

Consider a compact surface $S$ with a Riemannian metric with
constant negative curvature. Then the geodesic flow $\phi_t$
on the unit tangent bundle $M=T^1S$ of $S$ is an Anosov flow
whose strong stable and strong unstable foliations are $C^r$
smooth for all $r>1$; see e.g. Benoist, Foulon and
Labourie~\cite{BeFoLa92}.

In addition, both foliations are minimal: this type of
geodesic flow preserves a contact form, which coincides with
the Liouville measure on the unit tangent bundle, and it is
known to be ergodic with respect to this measure since the
work of Hopf (see e.g. \cite{Hopf71}). Therefore, the flow
is transitive and the entire phase space is non-wandering,
thus both strong stable and strong unstable foliations
(those tangent to $E^s$ and $E^u$ respectively) are minimal;
see Plante~\cite{Plante72}.

Hence we can apply our results to $f=\phi_1:M\to M$ which is
a strongly partially hyperbolic map and also preserves a
natural volume form that is a $u$-Gibbs state and the unique
physical measure for $f$ on $M$. We conclude that $f$ is
$C^2$ close to a $C^r$ stably ergodic (not necessarily
conservative) diffeomorphism with positive central Lyapunov
exponents Lebesgue almost everywhere.

\subsubsection{The time-1 map of the geodesic flow on
  symmetric Riemannian manifold of constant negative
  curvature}
\label{sec:time-1-map}

Anosov flows $\phi_t$ in any compact finite dimensional
Riemannian manifold having smooth (at least of class $C^3$)
strong stable or strong unstable foliations are essentially
$C^\infty$ conjugated to the geodesic flow over a locally
symmetric Riemannian manifold with constant negative
sectional curvature; see Benoist, Foulon and Labourie
\cite{BeFoLa92}. These flows preserve a smooth volume form
which is a contact form, so they are contact Anosov
flows. In addition, they are transitive by the classical
result of Hopf \cite{hopf39}, and do not admit sections; see
Godbillon~\cite[pp. 146-147]{Godb69}.  Hence, by the work of
Plante \cite{Plante72}, the strong unstable foliation is
minimal.

Hence we can apply our results to $f=\phi_1$ as in the
previous class of examples. We note that now the stable and
unstable directions are higher dimensional: if the dimension
of the manifold is $n$, the dimension of the unit tangent
bundle is $2n-1$, and then the dimension of the stable and
unstable invariant distributions equals $n-1$.

\subsection*{Acknowledgments}

Part of this work was done while V.A. was visiting
Universidad Cat\'olica del Norte, at Antofagasta, and
Pontificia Universidad Cat\'olica at Valpara\'{\i}so, Chile,
on several other occasions. V.A. wishes to thank these
institutions for their kind hospitality.




\section{Properties of $u$-Gibbs states
  and $cu$-Gibbs states}\label{sec:prelim-results}

Here we overview some properties of $u$- and $cu$-Gibbs states  used
in the detailed arguments of the previous section.

\subsection{Birkhoff regular points and ergodic
  decomposition}
\label{sec:birkhoff-regular-poi}

A point $z\in M$ is {\em Birkhoff regular} if  the Birkhoff averages
\begin{align}\label{eq:birkhoffneg}
 \vfi^-(z)&=\lim_{n\to\infty}\frac{1}{n}\sum_{k=0}^{n-1}
 \vfi(f^{-k}(z)),
 \\
 \vfi^+(z)&=\lim_{n\to\infty}\frac{1}{n}\sum_{k=0}^{n-1}
 \vfi(f^k(z));
\label{eq:birkhoffpos}
\end{align}
are defined and $\vfi^-(z)=\vfi^+(z)$ for every
$\vfi:M\to\RR$ continuous. The set of Birkhoff regular points of $f$ has full measure with respect to
any $f$-invariant measure $\mu$.

Given a point $x$ let us denote by $\mu_x$ the probability
measure given by the time average along the orbit of $x$
\begin{equation}\label{eq:ErgComp}
  \int\vfi\:d\mu_x=\lim_{n\to\infty}\frac{1}{n}\sum_{j=0}^{n-1}\vfi\big(
  f^j(x)\big)
\end{equation}
for every continuous $\vfi:M\to\RR$. According to the
Ergodic Decomposition Theorem $\mu_x$ is well defined and
ergodic for every $x$ in a set $\Sigma(f)\subseteq M$ that
has full measure with respect to any invariant measure
$\mu$; see e.g. Ma\~n\'e~\cite{Man87}. Moreover, for every
bounded measurable function $\vfi:M\to\RR$ we can write
\begin{equation}\label{eq:EDT}
 \int\vfi\:d\mu=\int\left(\int\vfi\:d\mu_x\right)d\mu(x).
\end{equation}
For every such $\vfi$ the integral $\int\vfi\:d\mu_x$
coincides with the time average $\mu$-almost everywhere, and
$x\in\supp \mu_x$ for $\mu$-almost all $x$.




\subsection{$u$-Gibbs states and their properties}\label{sec:existence-u-gibbs}

We assume from now on in this section that $f\::\:M\to M$ is
a $C^r$ partially hyperbolic diffeomorphism ($r\geq 2$)
having a splitting of the tangent bundle given by
$TM=E^s\oplus E^c\oplus E^u$. The aim of this subsection is
to present some usefull properties of $u$-Gibbs states.

An $f$-invariant probability measure $\mu$ is a {\em
  $u$-Gibbs state} if the conditional measures of $\mu$ with
respect to the partition into local strong unstable
manifolds are absolutely continuous with respect to Lebesgue
measure along the corresponding local unstable
manifolds. The following classical result shows that in this
setting there always exist $u$-Gibbs states.

\begin{proposition}\label{pr:uM1} \cite[by Pesin and Sinai]{PS82}
  Denote by $m_u=\dim E^ u$. If $D^u$ is an
  $m_u$-dimensional disk inside a strong-unstable leaf, and
  $\m_{D^u}$ denote the Lebesgue measure induced on $D^u$,
  then every accumulation point of the sequence of
  probability measures
$$\mu_n=\frac1n\sum_{j=0}^{n-1}f_*^j
\left(\frac{\m_{D^u}}{\m_{D^u}(D^u)}\right)$$ is a $u$-Gibbs
state with densities with respect to Lebesgue measure along
the strong-unstable leaves uniformly bounded away from zero
and infinity. In particular, the support of $\mu$ consists
of entire strong-unstable leaves.
\end{proposition}

\begin{remark}
  \label{rmk:uniform-density}
  As can be seen in Bonatti, D\'iaz and Viana \cite{BDV2004}
  the densities of $u$-Gibbs states with respect to Lebesgue
  measure along the strong unstable plaques depend only on
  $f$ through its derivatives and the curvature of the
  unstable manifolds. Consequently, the bounds on the
  densities of $u$-Gibbs states are also uniform for maps on
  a $C^2$ neighborhood of $f$.
\end{remark}

Next results present several very useful properties of
$u$-Gibbs states. First, the ergodic decomposition of
$u$-Gibbs states is formed by other $u$-Gibbs states.

\begin{proposition}\label{pr:uM2}{\cite[Lemma 11.13 and Corollary
    11.14]{BDV2004}}
  Ergodic components of any $u$-Gibbs state $\mu$ are
  $u$-Gibbs states whose densities are uniformly bounded
  away from zero and infinity. Conversely, a convex
  combination of $u$-Gibbs states is an $u$-Gibbs state. The
  support of any $u$-Gibbs state consists of entire
  strong-unstable leaves.
\end{proposition}

This allows us to assume without loss of generality in many
settings that $u$-Gibbs states are ergodic. The fact that
the support of any $u$-Gibbs state contains a full
strong-unstable leaf is very important when we assume that
the strong-unstable foliation is minimal, see
Section~\ref{sec:minimal-unstable-fol-1}.

Next we see that the ``basin of the family of all $u$-Gibbs
states'' is very big in the manifold.

\begin{proposition}\label{pr:uM6}{\cite[Theorem 11.16]{BDV2004}} 
  There exists $E\subseteq M$ intersecting every unstable
  disk on a full Lebesgue measure subset, such that for any
  $x\in E$, every accumulation point $\nu$
  of $\nu_{n,x}=(1/n)\sum_{j=0}^{n-1}\de_{f^j(x)}$ is
  a $u$-Gibbs state.
  \end{proposition}

  The following two propositions are consequences of
  Proposition~\ref{pr:uM6}.

\begin{proposition}\label{pr:uM3}{\cite[Section 11.2.3]{BDV2004}}
 If $\mu$ is an physical measure for $f$, then $\mu$ must
  be a $u$-Gibbs state.
\end{proposition}

\begin{proposition}\label{pr:uM4}\cite[by
  Dolgopyat]{Do2004-2} If $\mu$ is the unique $u$-Gibbs
  state for $f$, then $\mu$ is a physical measure for
  $f$. Moreover its basin $B(\mu)$ has full Lebesgue measure
  in $M$.
\end{proposition}

\subsection{$cu$-Gibbs states and their
  properties}\label{sec:existence-cu-gibbs}

The aim of this subsection is to present some usefull
results on $cu$-Gibbs states.

We recall that $E^{cu}:=E^c\oplus E^u$ and we denote
$m_{cu}:=\dim\ecu$. An invariant measure $\mu$ is a {\em
  $cu$-Gibbs state} if the $m_{cu}$ largest Lyapunov
exponents are positive $\mu$-almost everywhere and the
conditional measures of $\mu$ along the corresponding local
Pesin's center-unstable manifolds are $\mu$-almost
everywhere absolutely continuous with respect to Lebesgue
measure on these manifolds.

The notion of $cu$-Gibbs state was introduced by Alves,
Bonatti and Viana in \cite{ABV00} in a more general context;
they correspond to a non-uniform version of the $u$-Gibbs
states. In what follows, we present the properties of
$cu$-Gibbs states adapted to our setting. The interested
reader should consult~\cite{ABV00,BDV2004,vasquez2006} for
the properties of $cu$-Gibbs in more general contexts.

We first present a condition which guarantees the existence
of $cu$-Gibbs states. We say that a diffeomorphism $f$ has
{\em non-uniform expansion along the center-unstable
  direction} if there exists a constant $c_0>0$ such that
\begin{equation}\label{eq:mostly expanding}
\limsup_{n\to\infty}\frac{1}{n} \sum^{n-1}_{j=0}
\log\|Df^{-1}|E^{cu}(f^j(x))\|\leq -c_0<0.
\end{equation}
for all $x$ in a full Lebesgue measure subset of $M$.


%

\begin{proposition}\label{pr:ABV1}\cite[by Alves,
    Bonatti, Viana]{ABV00} 
  If $f$ is non-uniformly expanding along the center
  unstable direction, then there exist finitely many
  $cu$-Gibbs states $\mu_1,\dots,\mu_k$.  Moreover, they are
  physical measures a the union of their basins cover
  Lebesgue almost every point of the whole manifold $M$.
\end{proposition} 

If follows from Proposition~\ref{pr:uM6} that a $cu$-Gibbs
state is a $u$-Gibbs state. The converse is not true in
general, even if the $u$-Gibbs state is ergodic and it has
positive central Lyapunov exponents; see for instace the
example in \cite{vasquez2009}.

As explained along the statements of the main results, we
say that $f$ is {\em mostly expanding} if every $u$-Gibbs
state of $f$ has positive central Lyapunov exponents. As for
the dual notion of mostly contracting, the mostly expanding
diffeomorphisms possess similar properties which are
studied in \cite{AndVas2011}.

\begin{proposition}\label{mteo:MACA}{\cite[Theorem A]{AndVas2011}}
  The class of mostly expanding partially hyperbolic
  diffeomorphisms constitutes a $C^2$-open subset
  $\Diff$. Moreover, if $f$ is a mostly expanding partially
  hyperbolic diffeomorphism, then it is non uniformly
  expanding along the central direction so that, in
  particular, such $f$ has a finite number of physical
  measures whose basins together cover Lebesgue almost every
  point in $M$.
\end{proposition}

%

We can say more about the number of physical measures.

\begin{proposition}\label{pr:upper-semicont}{\cite[Theorem B]{AndVas2011}}
  If $f$ is a mostly expanding partially hyperbolic $C^r$
  diffeomorphism, $r\geq 2$, with a unique $cu$-Gibbs state,
  then every $g$ $C^r$-close to $f$ has a unique physical
  measure.
\end{proposition}

\subsection[Minimal unstable foliation and uniqueness of
$u$-Gibbs states]{Minimal unstable foliation, uniqueness of
  $u$-Gibbs states and stable ergodicity}
\label{sec:minimal-unstable-fol-1}

We now assume that $f$ is a partially hyperbolic $C^2$
diffeomorphism whose \emph{unstable foliation is 
  minimal}.

From the absolute continuity of the unstable foliation, we
see that the subset $E$ of $M$ given by
Proposition~\ref{pr:uM6} has full Lebesgue measure (volume)
in $M$.

So for Lebesgue almost every $x$ and for every given
continuous function $\vfi$ on $M$, we have that $\vfi^+(x)$
can have many different values, but each of them is given by
$\mu(\vfi)=\int \vfi\,d\mu$ where $\mu$ is some 
$u$-Gibbs state. 
 However, even if $f$ is a
strongly partially hyperbolic $C^2$ diffeomorphims with
\emph{simultaneously minimal stable and unstable
  foliations}, we cannot in general ensure that $\vfi^\pm$
is well defined nor that
$\vfi^+=\vfi^-$ Lebesgue almost everywhere.

These averages are well-defined Lebesgue almost everywhere and
 $\vfi^+=\vfi^-$ is true, for instance, if $f$
preserves Lebesgue measure.

\begin{problem}
  \label{prob:Birkhoffpastfuture}
  For a strongly partially hyperbolic $C^2$ diffeomorphims
  $f$ of a compact manifold $M$, if for any given continuous
  function $\vfi$ on $M$, both the forward $\vfi^+(x)$ and
  backward $\vfi^-(x)$ Birkhoff averages exist and coincide
  for Lebesgue almost every $x\in M$, then Lebesgue measure
  is invariant.
\end{problem}

Moreover, since every strong-unstable leaf is dense by assumption,
and the support of every $u$-Gibbs state contains some full
strong-unstable leaf, we deduce from Proposition~\ref{pr:uM6} that,
for each $x\in E$, every accumulation measure $\mu_x$ given by
\eqref{eq:ErgComp} has full support. In particular
$\mu_x(U)>0$ for every open subset $U$. 

Altogether this ensures that each $x\in E$ has positive
frequency of visits to any open subset of the manifold. In
particular, \emph{Lebesgue almost every positive orbit is
  dense}. 

We will construct a perturbed map $g$ $C^2$ close to $f$
whose strong-unstable foliation coincides with the
strong-unstable foliation of $f$, so \emph{these properties
  persist for all maps $g$ obtained from $f$ according to
  our perturbation scheme}, to be presented in the next
section.

A similar argument shows that \emph{physical $cu$-Gibbs
  states are unique}, if they exist.

\begin{lemma}
  \label{le:unque-cuGibbs}
  Let $f:M\to M$ be a (strongly) partially hyperbolic $C^2$
  diffeomorphism whose strong-unstable foliation $\cF^u$ is
  minimal. Is $f$ is mostly expanding, then there can exist at most one $cu$-Gibbs state
  for $f$ which is also a physical measure and it basin
  cover Lebesgue almost every point of the whole manifold
  $M$.
\end{lemma}

\proof Indeed, an ergodic $cu$-Gibbs state $\mu$ is a physical
measure whose basin $B(\mu)$ contains some open neighborhood
$U$ of any center-unstable  leaf $\xi$ supporting $\mu$. These
neighborhoods are given by the family of strong-stable
leaves through points of $\xi$. Therefore any other ergodic
$cu$-Gibbs state $\nu$ has some  center-unstable leaf $\zeta$ is its support
which crosses $U$ (Note that in  this case, the center unstable leaf $\zeta$ is foliated by strong unstable leafs and they are dense in $M$). 
The absolute continuity of the
strong-stable foliation ensures that the same argument as
above implies $\mu=\nu$.

The basin of $\mu$ has full Lebesgue measure in the
ambient space as a direct consequence of
Proposition~\ref{pr:uM6}. For otherwise it would be possible
to find an invariant subset of positive Lebesgue measure
where \eqref{eq:mostly expanding} holds and then, using the
proof of Proposition~\ref{pr:ABV1}, we would be able to
construct another $cu$-Gibbs state which would be a new
physical measure.\endproof

%
%
%
%

Combining the previous results we have proved the following.

\begin{proposition}
  \label{pr:stbl-erg-most-contr}
  Let $f$ be a mostly expanding strongly partially
  hyperbolic diffeomorphism, with one-dimensional
  subbundles, having a unique $cu$-Gibbs state.

  Then there exists a $C^r$ neighborhood $\U$ of $f$,
  $r\ge2$, such that each $g\in\U$ is mostly expanding and
  admits an unique $cu$-Gibbs state $\mu_g$ which is a
  physical measure with full basin. In particular $f$ is
  $C^r$-stably ergodic.
\end{proposition}

\section{Perturbing the central subbundle}
\label{sec:pert-centr-lyap}

Here we start the proof of Theorem~\ref{mthm:removing},
providing the details of the construction of the perturbed
diffeomorphism.

\subsection{Adapted norms for partially hyperbolic diffeomorphisms}
\label{sec:adapted-norms-partia}

We note first that for a flow $X_t$ without equilibria on a
compact Riemannian manifold $M$ with an induced norm
$\|\cdot\|$, we can define a new norm
$|u|_x:=\|u\|/\|X(x)\|$ which satisfies
\begin{align*}
  |DX_t\cdot X(x)|_{X^t(x)}
  =
  \frac{\|DX_t\cdot X(x)\|}{\|X(X_t(x))\|}
  =1, \quad x\in M,  t\in\RR.
\end{align*}
We then use the results from Gourmelon~\cite{Goum07}, which
provide adapted metrics on our setting. That is, we may
assume without loss of generality that our partial
hyperbolic map $f=X_1$, where $(X_t)_{t\in\RR}$ is now a
$C^2$ Anosov flow, satisfies \eqref{PH1}, \eqref{PH2} and
\eqref{PH3} with $C=1$ and, moreover, that the norm is
induced by a Riemannian metric on $M$ such that at each
$x\in M$ the directions along the different subbundles are
mutually orthogonal.

Since $f=X_1$ is the time-$1$ map of an Anosov flow, we have
$$\|Df(x)v\|=1
\quad\text{for every $v\in E^c_f(x)$, with $\|v\|=1$.}$$

\subsection{The choice of coordinates}
\label{sec:dynam-coher-dynam}

The assumption that $f=X_1$ is the time-$1$ map of a $C^2$
Anosov flow ensures that $E^c_f$ is an integrable subbundle,
its integral manifolds $\cF^c$ are the orbits of the flow
$(X_t)_{t\in\RR}$. Moreover the strong-stable $\cF^{s}$ and
strong-unstable $\cF^{u}$ foliations are smooth along the
central leaves since
$\cF^{\star}(X_t(x))=X_t(\cF^{\star}(x))$ for $x\in M,
t\in\RR$ and $\star=s,u$. This means that we have
\emph{dynamical coherence}: for each $x\in M$ the sets
$\cF^{cu}(x)=\cup_{t\in\RR}X_t(\cF^{u}(x))$ and
$\cF^{cs}(x)=\cup_{t\in\RR}X_t(\cF^{s}(x))$ are $C^2$
immersed submanifolds of $M$ tangent to $E^{cu}:=E^c\oplus
E^u$ and $E^{cs}:=E^s\oplus E^c$ respectively.

Using the assumption of $C^2$ smoothness of the
strong-unstable foliation $\cF^{u}$ together with the
dynamical coherence we see that the \emph{center-unstable
  foliation $\cF^{cu}$ is also $C^2$ smooth}.  Hence we can
define for any given $q_0\in M$ a $C^2$ local
parametrization $\psi:D\to M$ of $M$ as follows.

We start by choosing an embedded manifold $\Sigma$ in a
neighborhood of $q_0$, containing $q_0$, such that $\Sigma$
is transverse to $E^u_f$ at all points $q\in\Sigma$: we
simply take
$\Sigma:=\cup_{-\gamma<t<\gamma}X_t(W^{s}_{\gamma}(q_0))$
where $W^s_\gamma(q_0)$ is the connected component of
$\cF^s(x)\cap B_\gamma(q_0)$ containing $q_0$. This is a
local center-stable leaf through $q_0$.

We write $W^{cs}_\gamma: \widetilde B_1(0)\to\Sigma$ for a
$C^2$ parametrization of $\Sigma$, where $\widetilde B_1(0)$
is the unit ball in the Euclidean space $\RR^{s+c}$, with
$s:=\dim E^s_f$ and $c:=\dim E^c_f=1$ and
$W^{cs}_\gamma(0)=q_0$. We can assume that
\begin{align*}
  DW^{cs}(0)(\RR^s\times\{0^c\})=E^s(q_0)
  \qand
  DW^{cs}(0)(\{0^s\}\times \RR^c)=E^c(q_0).
\end{align*}
In this way $\Sigma$ is tangent to the center-stable
direction at all points, that is,
$T_{q}\Sigma=E^{cs}(q):=E^s_f(q)\oplus E^c_f(q),
q\in\Sigma$.

The Stable/Unstable Manifold Theorem ensures that for every
$q\in\Sigma$ there exists a $C^2$ embedding
$W^u_\gamma(q):\widehat B_1(0)\to M$ where $\widehat B_1(0)$ is
the unit ball in the Euclidean space $\RR^{u}$ with $u:=\dim
E^u_f$; $W^u_\gamma(q)(\widehat B_1(0))=\cF^u(q)\cap
B_\gamma(q)$ and $T_q W^u_\gamma(q)=E^u_f(q)$; see e.g. \cite{Sh87}.

The $C^2$ smoothness of $\cF^u$ ensures that the following
is a $C^2$ map
\begin{align*}
  \psi: \widetilde B_1(0)\times\widehat B_1(0)\to M,
  \quad
  (w,z)\mapsto W^u_\gamma\big(W^{cu}_\gamma(w)\big)(z)
\end{align*}
and, since $D\psi(w,0):\RR^{s+c+u}\to T_qM$ is an
isomorphism for all $w\in \widetilde B_1(0)$ and $\psi$ maps
$\widetilde B_1(0)\times\{0\}$ diffeomorphically onto
$\Sigma=W^{cu}_\gamma(\widetilde B_1(0))$, then $\psi$ is a
$C^2$ diffeomorphism on a neighborhood of $\widetilde
B_1(0)\times\{0\}$. Thus, setting $\gamma>0$ smaller if
needed, we can assume without loss that $\psi$ is a $C^2$
diffeomorphism between $D=\widetilde B_1(0)\times\widehat
B_1(0)$ and its image in $M$.  

This is the parametrization we need to define our
perturbation.  We note that
\begin{align*}
  D\psi(w,0)(\RR^{s+c}\times\{0^u\})
  &=
  E^s_f(\psi(w,0))\oplus E^c_f(\psi(w,0)) 
  \quad\text{and}
  \\
  D\psi(w,z)(\{0^{s+c}\}\times\RR^u)
  &=
  E^u_f(\psi(w,z))
  \quad\text{and also}
  \\
  D\psi(w,z)(\{0^{s}\}\times\RR\times\{0^u\})
  &=
  E^c_f(\psi(w,z)),
\end{align*}
for all $(w,z)\in D$. The last property is a consequence of
the assumption that $f=X_1$, since for all $(w,z)\in D$ and
each small $|t|$ there exists $z^\prime$ such that
$\psi(X_t(q),z)=X_t(\psi(q,z^\prime))$.

So this provides a $C^2$ coordinate system around each point
of $M$ such that $(w,z)\mapsto\psi(w,z)$, for $(w,z)\in D$,
is contained in the local strong-unstable manifold of the
point $\psi(w,0)$.

\subsection{Construction of the local perturbation}
\label{sec:constr-local-perturb}

Let us fix $q_0$ a recurrent non-periodic point of $f$. We
note that since we assume $f$ has a minimal unstable
foliation, then we have plenty of points with dense
orbit. We fix a neighborhood $U$ of $q_0$ in $M$ such that
$f(U)\cap U=\emptyset$ and a parametrization $\psi:D\to U$,
where $\psi$ is the coordinate system constructed in the
previous subsection.
By appropriately rescaling the basis vectors in $\RR^{d}$, we
can further assume without loss of generality that
$\|D\psi(0)e_i\|=1, i=1,\dots,d$, for the canonical basis
$\{e_i\}_{i=1}^{d}$ of $\RR^d$, where $d:=s+c+u=\dim M$.

\subsubsection{The choice of the bump function}
\label{sec:choice-bump-functi}

\begin{figure}[htpb]
  \centering
  \includegraphics[width=5cm]{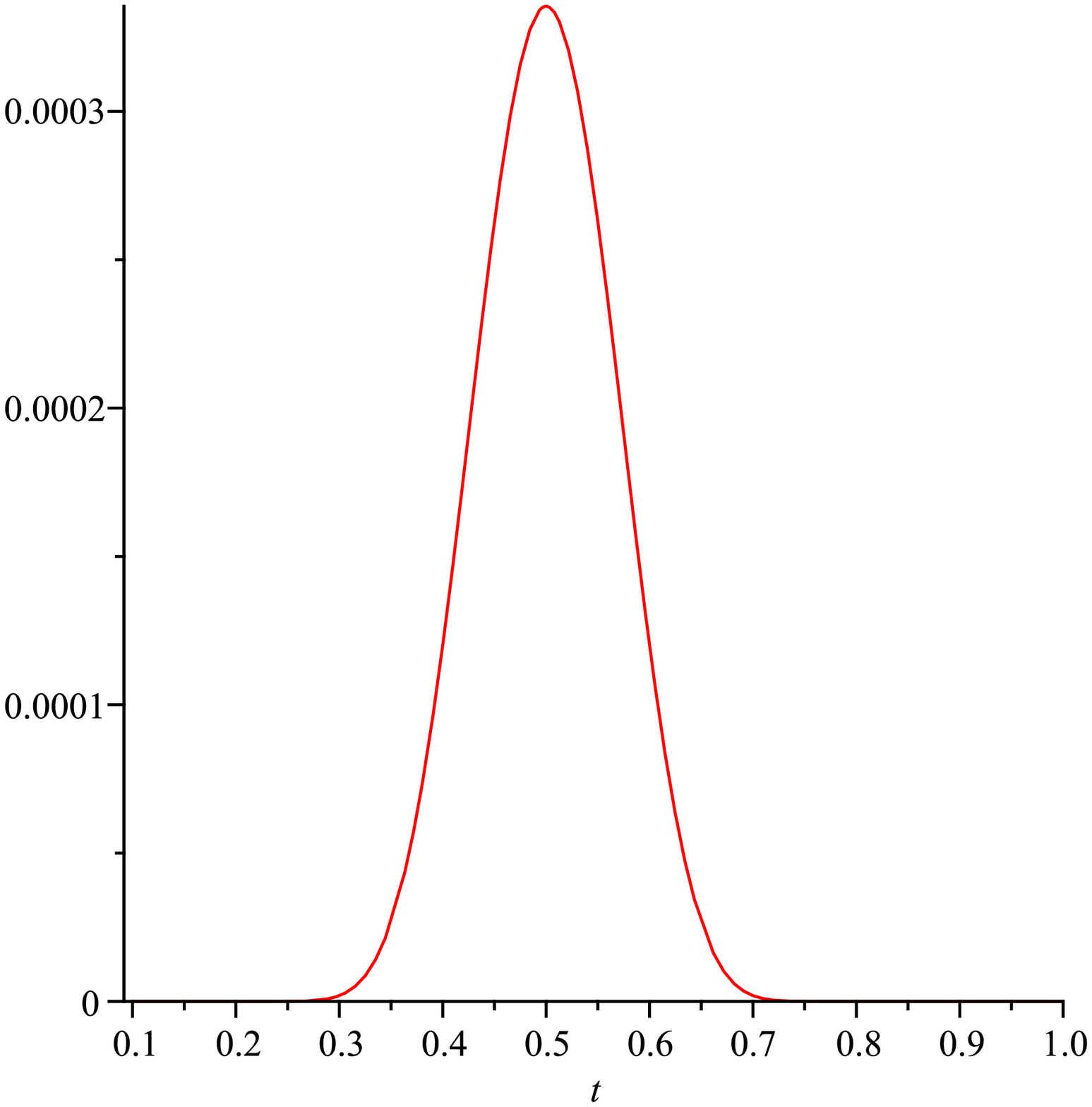}
  \includegraphics[width=5cm]{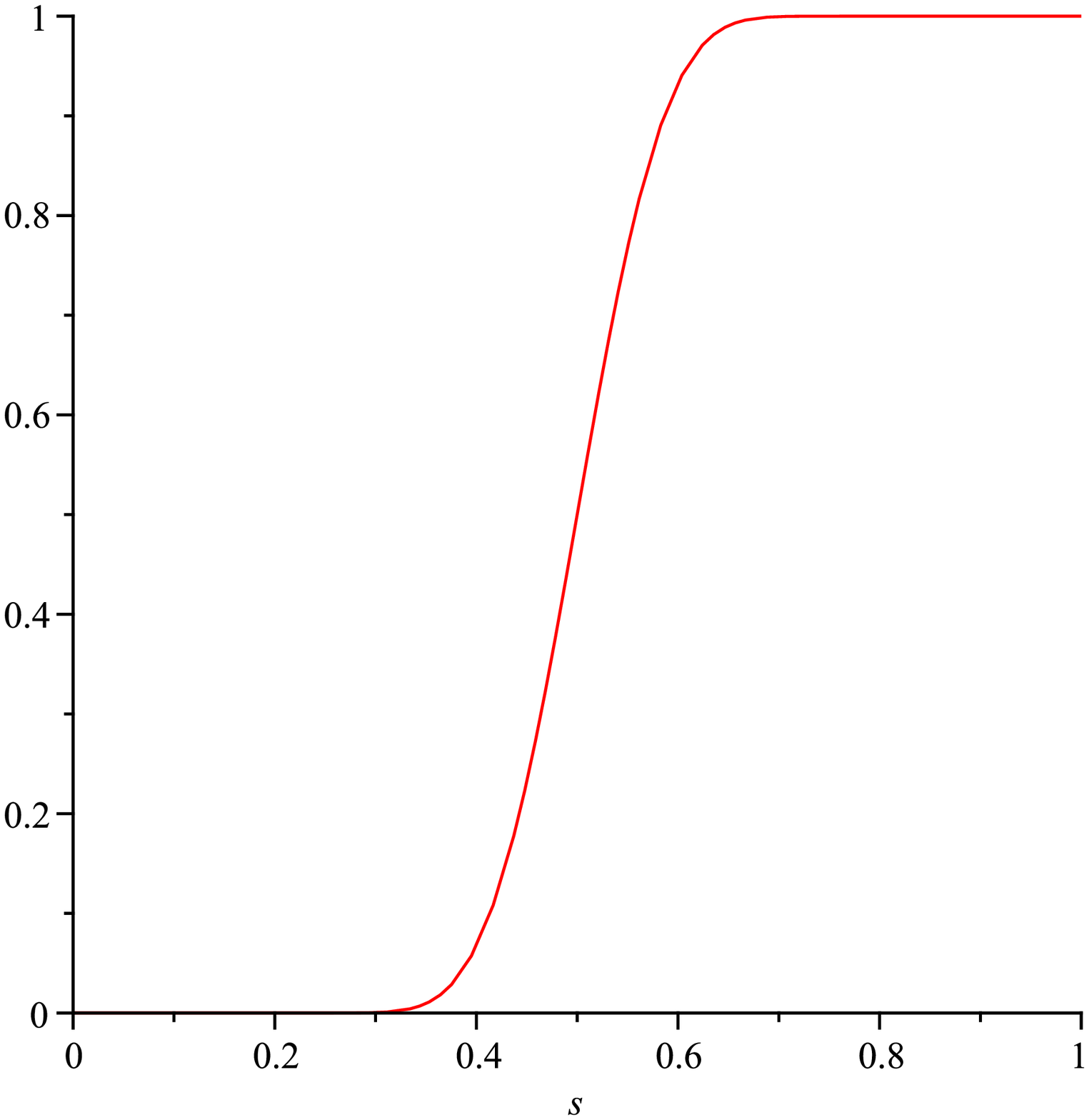}
  \\
  \includegraphics[width=5cm]{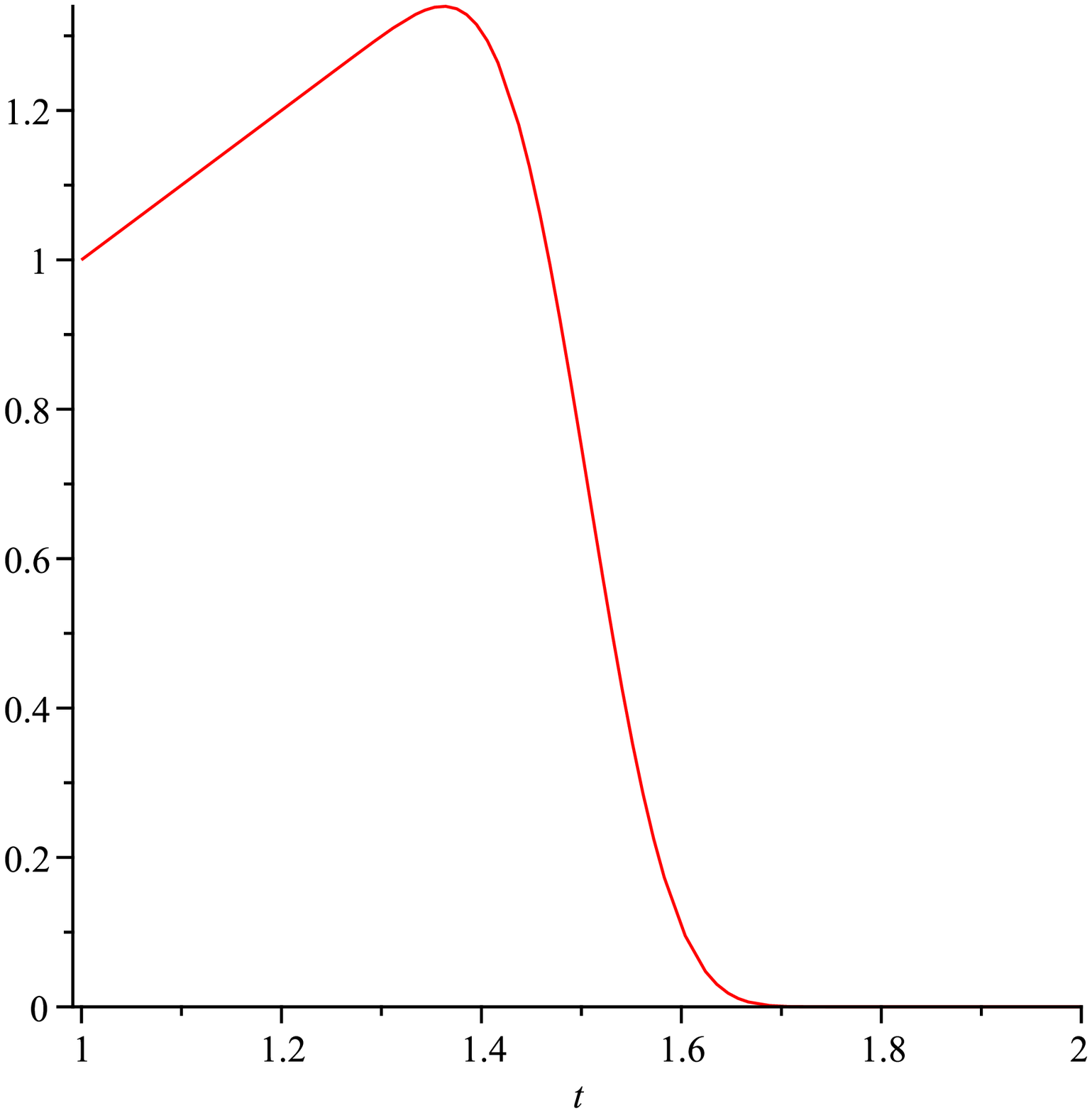}
  \includegraphics[width=5cm]{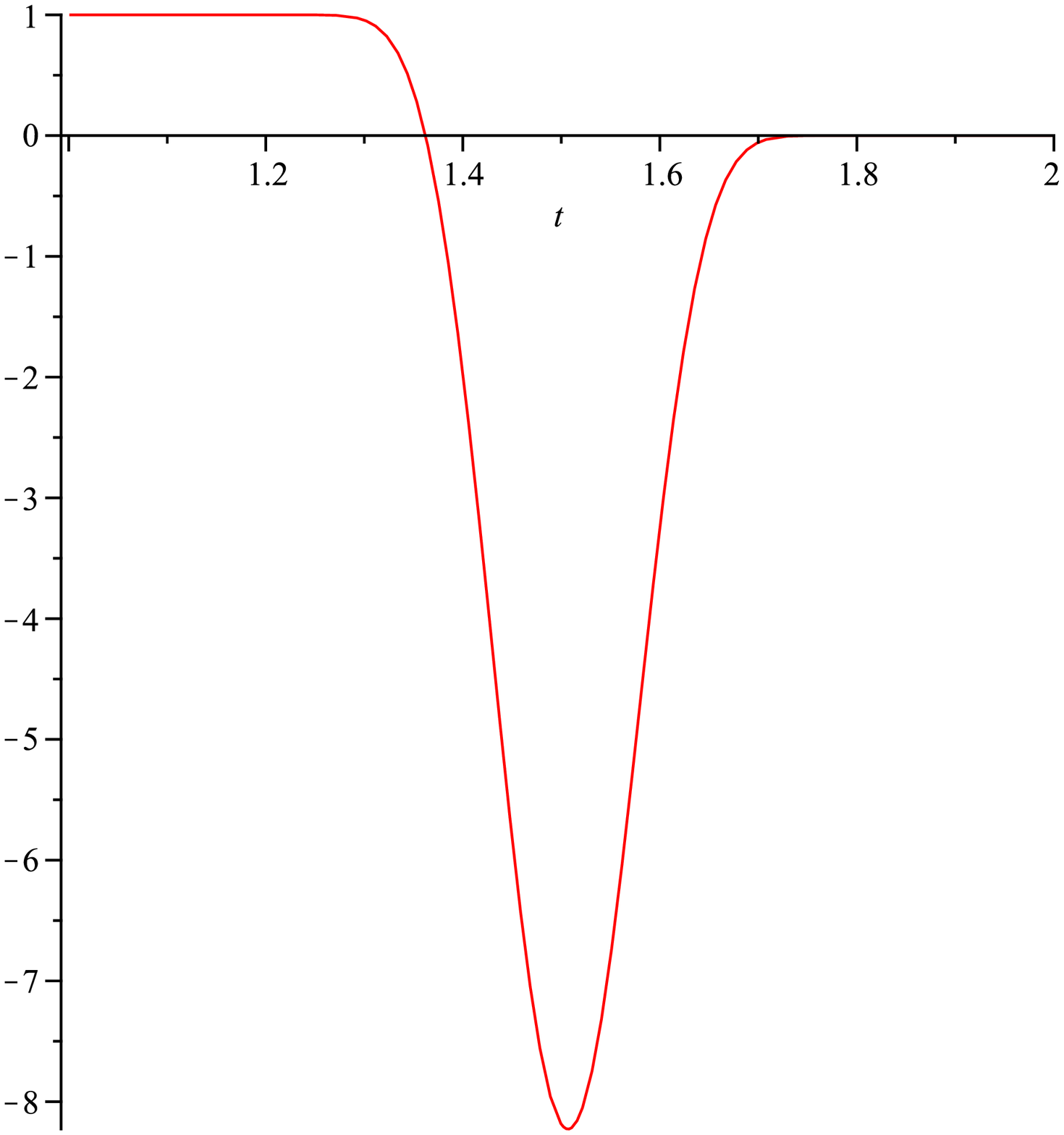}
  \caption{At the upper left we have $\eta_1$ and at the
    upper right $\eta_2$ on the interval $[0,1]$. At the
    lower left we have $\zeta(s)=s\phi(s)$ and at the lower
    right $\zeta^\prime$ on the interval $[1,2]$.}
  \label{fig:bump}
\end{figure}

Let $\phi:\RR\to[0,1]$ be a bump function with the following
properties: 
\begin{itemize}
\item $\phi\equiv0$ on $\RR\setminus(-2,2)$ and $\phi\equiv 1$
  on $[-1,1]$;
\item $(s\phi(s))^{'}\neq0$ for all $s\in(-2,2)$ except at
  the points $\{\pm s_0\}$ for a value $s_0\in(1,2)$.
\end{itemize}
This is easy to build: we start with
\begin{align*}
  \eta_0(s):=
  \begin{cases}
    e^{-1/t^2} & \text{if } t>0 \\
    0 &  \text{if } t\ge0
  \end{cases}
  \qand
  \eta_1(s):=\eta_0(s)\eta_0(1-s)
\end{align*}
and then consider $\eta_2(s):=c^{-1}\int_{-\infty}^s
\eta_1(t)\,dt$ where
$c:=\int_{-\infty}^\infty\eta_1(t)\,dt$. This function
$\eta_2$ is of class $C^\infty$ and satisfies $\eta_2\mid(-\infty,0]\equiv0$,
$\eta_2\mid[1,+\infty)\equiv1$ and $\eta_2\mid(0,1)>0$.
We set $\phi(s):=\eta_2(s+2)\eta_2(2-s)$. See
Figure~\ref{fig:bump}.

For $s\in(-1,1)$ we have $s\phi(s)=s$ and there is no zero
of the derivative of $\zeta(s):=s\phi(s)$. But
$\zeta(0)=\zeta(2)=0$, thus there exists a zero of
$\zeta^\prime$ in $(1,2)$.  Since $\zeta(-s)=\zeta(s)$ we
may consider only $s\in(1,2)$ and clearly see in
Figure~\ref{fig:bump} that there exists a unique zero of
$\zeta^\prime$ in this interval, as we wanted.

\subsubsection{The choice of the perturbation}
\label{sec:choice-perturb}

To define the perturbation we set the
multi-indexes $x=(x_1,\dots,x_s)$ and $z=(z_1,\dots,z_u)$.
We then define
\begin{align*}
  \Phi_\epsilon(x,y,z):=
  \left[\prod_{i=1}^s\phi(\frac{x_i}{\epsilon})\right]\phi(\frac{y}{\epsilon})
  \left[\prod_{j=1}^u\phi(\frac{z_j}{\epsilon})\right]
\end{align*}
together with the following diffeomorphism for small $t>0$
(to be bounded above in the following arguments) and
$0<\epsilon<1/4$
\begin{align}\label{eq:def-h}
  h=h_{t,\epsilon}:D\to D, \quad (x,y,z)\mapsto
  (x,y,z+ty\Phi_\epsilon(x,y,z)e_{s+c+1}),
\end{align}
where $s+c+1$ is the first coordinate along the unstable
direction in the parametrization $\psi$.  We have, with
respect to the canonical basis on $\RR^d$ and writing $I_k$
for the identity on $\RR^k$ for $k\in\ZZ^+$
\begin{align}\label{eq:Dh}
  Dh_{t,\epsilon} =
  \begin{bmatrix}
    I_s & 0 & 0 \\
    0 & 1 & 0 \\
    t y D_x\Phi_\epsilon \cdot e_{s+c+1} &
    t(\Phi_\epsilon+y\partial_y\Phi_\epsilon) \cdot e_{s+c+1} & 
    I_u+t y D_z\Phi_\epsilon \cdot e_{s+c+1}
  \end{bmatrix},
\end{align}
where
$D_x\Phi_\epsilon(x,y,z): \RR^s\to\RR$ and
$D_z\Phi_\epsilon(x,y,z):\RR^u\to\RR$. Clearly $h$ is the identity on
$D\setminus B(0,2\epsilon)$.  For $(x,y,z)\in
B(0,2\epsilon)$ and some constant $C>0$ (a bound on
$\sup|D\phi|$) we have for $\star=$ any of the variables in $x,y,z$
\begin{align*}
  |ty \partial_\star \Phi_\epsilon|\le t\cdot 2\epsilon (C/\epsilon)=2Ct
  \quad\text{and}\quad
  |t(\Phi_\epsilon+y\partial_y\Phi_\epsilon)|\le t(1+2C),
\end{align*}
and by the definition of $\phi$ we get
$|t(\Phi_\epsilon+y\partial_y\Phi_\epsilon)|\neq0$ except at
two values $y=\pm\epsilon s_0$ with
$|y|\in(\epsilon,2\epsilon)$. We note that the bottom right
side block of the matrix in \eqref{eq:Dh} has a non-zero
determinant
\begin{align}\label{eq:det-dh}
  |\det Dh_{t,\epsilon}|= | 1+ty\partial_{z_1}\Phi_\epsilon|\ge
  1-2Ct \ge \frac12 \quad\text{for $0<t<1/4C$.}
\end{align}
Moreover, since $\|h-Id\|=|ty\Phi_\epsilon|$, we see that
for small $\epsilon>0$ we get
$\|h_{t,\epsilon}-Id\|_{C^1}\le (1+2C)t$.
In addition, for $0<t<1/4C$ we have
\begin{align*}
  h(x,y,z)=h(\bar x,\bar y,\bar z)
  \implies
  x=\bar x, y=\bar y , z_j=\bar z_j
  \quad\text{for}\quad
  j=s+c+2,\dots,d
  \end{align*}
and
\begin{align*}
  z_{s+c+1}-\bar z_{s+c+1} = t y [ \Phi_\epsilon(x,y,\bar
  z)-\Phi_\epsilon(x,y,z)].
\end{align*}
From the definition of $\Phi_\epsilon$ we see that
$\Phi_\epsilon(x,y,\bar z)=\Phi_\epsilon(x,y,z)=0$ for
$|y|>2\epsilon$; while for $|y|\le2\epsilon$
\begin{align*}
  |z_{s+c+1}-\bar z_{s+c+1}|
  &\le 
  t|y| 
  \left| 
    \phi\Big(\frac{\bar z_{s+c+1}}\epsilon\Big)-
    \phi\left(\frac{z_{s+c+1}}{\epsilon}\right) \right|
  \\
  &\le
  t\frac{|y|}{\epsilon}(\sup|D\phi|)|z_{s+c+1}-\bar z_{s+c+1}|
  \le
  2t C |z_{s+c+1}-\bar z_{s+c+1}|
\end{align*}
so from the assumption $0<t<1/4C$ we conclude that $z=\bar
z$ in all cases. Under this condition the map $h$ is
injective and thus, by \eqref{eq:det-dh}, a diffeomorphism
onto its image.

\begin{remark}\label{rmk:higher-diff}
  We can easily bound the second partial derivatives as
  $|\partial^2_{\star\dag}\Phi_\epsilon|\le C\epsilon^{-2}$
  for $\star,\dag=$ the variables in $x,y,z$, so that
  $\|h_{t,\epsilon}-Id\|_{C^2}\le Ct\epsilon^{-1}$ for a
  constant $C>0$ depending on the bump function $\phi$ but
  independent of $\epsilon,t$. Moreover, for $r>2$ we have
  $|\partial^k_m\Phi_\epsilon|\le C_k\epsilon^{-k}$ for
  every multi-index
  $m=(m_1,\dots,m_k)\in\{x_1,\dots,x_s,y,z_1,\dots,z_u\}^k$
  and $2\le k\le r$, where $C_k$ depends on $\phi$
  only. This will be essential to estimate the $C^r$
  distance of the perturbed map away from $f$.
\end{remark}

We finally define the perturbed map $g$ as
\begin{align*}
  g(x)=
  \begin{cases}
    f(q) & \text{if  } q\in M\setminus U\\
    (f\circ H)(q) & \text{if  } q\in U
  \end{cases},
\end{align*}
where we write $H=\psi\circ h\circ \psi^{-1}$ from now on.

\begin{remark}\label{rmk:same-return}
  We note that $g(U)=f(U)$ by the choice of $U$ since
  $H(U)=U$. Thus we see that $g(U)\cap U = \emptyset$ and
  that the minimum $n>0$ so that $g^n(U)\cap U\neq\emptyset$
  is at least $2$ and \emph{depends only on $f$}, because
  $g=f$ outside of $U$.
\end{remark}

We observe that $Dh(0)(u,v,w)=(u,v,w+tv\cdot e_{s+c+1})$ so
at $q_0$ we have $H(q_0)=q_0$ and the image $Dg(q_0)\cdot
E_f^c(q_0)$ is the graph of a non-zero injective linear map
$L_{g(q_0)}:E^c_f(g(q_0))\to E^u_f(g(q_0))$.

In the absence of dynamical coordinates, the splitting
$E^s_f\oplus E^c_f\oplus E^u_f$ in general depends not more
than H\"older continuously on the base point. Hence it is
not possible in general to find a smooth coordinate change
that sends the $Df$-invariant direction onto the coordinate
axis everywhere in a neighborhood of $q_0$. But the choice
of $\psi$ through dynamical coherence ensures that for $q\in
V=V_\epsilon:=\psi(B(0,2\epsilon))$ the unstable and
center-unstable directions are preserved by $DH(q)$.  Hence
we can ensure that \emph{there is a non-zero injective
  linear map $L_{g(q)}:E^c_f(g(q))\to E^u_f(g(q))$ such that
  $Dg(q)\cdot E^c_f(q)$ is the graph of $L_{g(q)}$}.

In particular, if $\pi^u_{g(q)}:T_{g(q)}M\to
E^u(g(q))$ is the projection into $E^u_f(g(q))$ parallel to $E^s_f(g(q))\oplus E^c_f(g(q))$ then it is
non-zero, i.e., $\pi^u_{g(q)}\circ L_{g(q)} \not\equiv 0$;
but the projection $\pi^s_{g(q)}:T_{g(q)}M\to E^u(g(q))$
parallel to $E^c_f(g(q))\oplus E^u_f(g(q))$ is zero, i.e.,
$\pi^s_{g(q)}\circ L_{g(q)}\equiv 0$.
This means that \emph{the image of the central vectors under
  $Dg$ has a non-zero unstable component along the original
  splitting but can only have a zero stable component.}
We remark that
\begin{itemize}
\item for $q\in M\setminus V$ the old invariant directions
  are preserved by $Dg$, i.e.  $Dg(q)\cdot
  E^*_f(q)=E^*_f(g(q))=E^*_f(f(q))$ since $Dg(q)=Df(q)$, for
  each $*=s,c,u$;
\item for $q\in V$ both the direction $E^u_f$ and $E^{cu}_f$
  are preserved by $Dg$, i.e., $Dg(q)\cdot
  E^u_f(q)=E^u_f(g(q))=E^u_f(f(H(q)))$ and $Dg(q)\cdot
  E^{cu}_f(q)=E^{cu}_f(g(q))=E^{cu}_f(f(H(q)))$.
\end{itemize}
Consequently, \emph{the unstable and center-unstable
  subbundles of $f$ remain as unstable and center-unstable
  subbundles for $g$: $E^u_g=E^u_f$ and
  $E^{cu}_f=E^{cu}_g$. Consequently the strong unstable and
  center-unstable foliations of $f$ and $g$ coincide since
  these foliations are uniquely integrable in our setting
  and, in particular, $\F^u_g$ is minimal}.

\subsection{The perturbed map is $C^2$ close}
\label{sec:perturb-map-c2}

We can perform all the previous constructions with a family
$h_{t,\epsilon}$ with $\epsilon$ going to zero and a
function $t=t(\epsilon)$ which also goes to zero, but
essentially arbitrary; see the next sections.  This freedom
of choice for $t(\epsilon)$ enables us to control the
distance of $g$ to $f$ in the $C^2$ topology.

Indeed, if the strong unstable foliation $\cF^u_f$ of $f$ is
of class $C^r$, for some $r>2$, then we can build $H$ of
class $C^r$ and, from Remark~\ref{rmk:higher-diff}, we have
$\|h_{t,\epsilon}-Id\|_{C^r}<C_rt\epsilon^{-r}$. So we just
have to choose the appropriate function $t(\epsilon)$.  In
the present scenario, we have a strong unstable foliation
$\cF^u_f$ of $f$ is of class $C^2$, and we may choose
$t(\epsilon)=\min\{\epsilon^3,1/(4C)\}$.  With this choice
of $t=t(\epsilon)$ we have
\begin{align*}
  \|h_{t,\epsilon}-Id\|_{C^2}
  \le 
  C_2 \frac{\epsilon^3}{\epsilon^{2}}\xrightarrow[\epsilon\to0]{}0.
\end{align*}
So in what follows we assume that we have performed the
perturbation described in the previous subsections with the
aid of the family of functions
$(h_{t(\epsilon),\epsilon})_{\epsilon\ge0}$ where
$t(\epsilon)$ was defined above.


\section{The perturbed central subbundle}
\label{sec:perturb-central-subb}

Here we prove the following lemma.

\begin{lemma}\label{le:Ecg-away-Ecf}
  For all small enough $0<\epsilon<1/4$ and $t=t(\epsilon)$,
  there exists a subset $\tilde V$ of $V$ such that:
  \begin{enumerate}
\item for each $u$-Gibbs state $\mu$ of $g$, $\mu(\tilde V) >0$;
\item for every $q\in \tilde V$,  $E^c_g(q)\ne E^c_f(q)$ .
  \end{enumerate} 
  In particular, if $q\in \tilde V$, $E^c_g(q)$ can be
  written as the graph of a nonzero linear map
  $G_q:E^c_f(q)\to E^u_f(q)$.
\end{lemma}

\proof Fix $q\in M$. For $v\in T_q M$, denote
$v^u:=\pi^u_q(v)$; 
$v^c:=\pi^s_q(v)$ 
and $v^s:=\pi^s_q(v)$ 
as defined in Section~\ref{sec:choice-perturb}.

For $v=v^c+v^u\in T_qM$, the slope of $v$ is defined as
$$s_c(v)=\frac{\|v^u\|}{\|v^c\|}.$$

We want to estimate the the slope of $Dg^n(q)v$ for, $n\geq
0$, $q\in M$ and $v=v^c+v^u\in E^c_f(q)\oplus E^u_f(q)$.
Note that if $q\in V$, then $g=f\circ H$, and so
\begin{eqnarray}\label{eq:slopeinV}
  s_c(Dg(q)v)&=&\frac{\|[Dg(q)v]^u\|}{\|[Dg(q)v]^c\|}=\frac{\|Df(H(q))[DH(q)v]^u\|}{\|Df(H(q))[DH(q)v]^c\|} \nonumber\\ \nonumber \\
  &=& \frac{\|Df(H(q))[DH(q)v]^u\|}{\|[DH(q)v]^u\|}\frac{\|[DH(q)v]^u\|}{\|[DH(q)v]^c\|}\frac{\|[DH(q)v]^c\|}{\|Df(H(q))[DH(q)v]^c\|}\nonumber\\ \nonumber\\
  &=& \frac{\|Df(H(q))[DH(q)v]^u\|}{\|[DH(q)v]^u\|}s_c(DH(q)v).
\end{eqnarray}
The last equality is obtained from $\|Df(q)v^c\|=\|v^c\|$
for every $q\in M$ and every $v^c\in E^c_ f(q)$ (the central
direction of $f$ is the flow direction).

On the other hand, if $q\in M\setminus V$, since $f=g$ then  
\begin{eqnarray}\label{eq:slopeoutV}
  s_c(Dg(q)v)&=&\frac{\|[Dg(q)v]^u\|}{\|[Dg(q)v]^c\|}=\frac{\|Df(q)v^u\|}{\|Df(q)v^c\|}\nonumber\\ \nonumber\\
  &=& \frac{\|Df(q)v^u\|}{\|v^u\|}\frac{\|v^u\|}{\|v^c\|}\frac{\|v^c\|}{\|Df(q)v^c\|}= \frac{\|Df(q)v^u\|}{\|v^u\|}s_c(v).
\end{eqnarray}

If $q\in V$, $g(q),\dots, g^{n}(q)\notin V$ and
$g^{n+1}(q)\in V$, $n\geq 1$, combining \eqref{eq:slopeinV}
and \eqref{eq:slopeoutV} we obtain,
\begin{eqnarray}\label{eq:accslope}
  s_c(Dg^{n+1}(q)v)&=&\frac{\|Df^n(g(q))[Dg(q)v]^u\|}{\|[Dg(q)v]^u\|}s_c(Dg(q)v)\nonumber \\ \nonumber\\
  &=&\frac{\|Df^n(g(q))[Dg(q)v]^u\|}{\|[Dg(q)v]^u\|}\frac{\|Df(H(q))[DH(q)v]^u\|}{\|[DH(q)v]^u\|}s_c(DH(q)v)\nonumber \\ \nonumber\\
  &=&\frac{\|Df^{n+1}(H(q))[DH(q)v]^u\|}{\|[DH(q)v]^u\|}s_c(DH(q)v).
\end{eqnarray}
In particular, we obtain the following bounds for the slope
\begin{equation}\label{eq:slopebounds1} {\rm
    m}(Df^{n+1}|E^u_f)s_c(DH(q)v) \leq s_c(Dg^{n+1}(q)v)\leq
  \|Df^{n+1}|E^u_f\|s_c(DH(q)v),
\end{equation}
and from \eqref{PH3} we can write
\begin{equation}\label{eq:slopebounds2}
  \lambda_3^{n+1}s_c(DH(q)v) \leq s_c(Dg^{n+1}(q)v)\leq \mu_3^{n+1}s_c(DH(q)v).
\end{equation}

For $q=q_0$, we have $H(q_0)=q_0$ and
$DH(q_0)(0,v^c,v^u)=(0,v_c,tv_c+v^u)$. If $v=v^c+v^u\in
E^{cu}_f(q_0)$, then
\begin{equation}\label{slopebounds3}
  s_c(DH(q_0)v)=\frac{\|tv^c+v^u\|}{\|v^c\|}\geq|t|-s_c(v).
\end{equation}
In particular, if $s_c(v)<\frac{|t|}4$,
then $$s_c(DH(q_0)v)>\frac{3|t|}4>\frac{|t|}4>s_c(v).$$ By
the continuous depence of the splitting with respect to the
point $q$ for $f$, there exist $\eta>0$ and $a>b>0$ such
that, for any $q\in B(q_0,\eta)$ and $v\in E^{cu}_f(q)$,
such that $s_c(v)<b$, we have
\begin{equation}\label{eq:slopebounds4}s_c(DH(q)v)>a>b>s_c(v).\end{equation}

\begin{figure}[h]
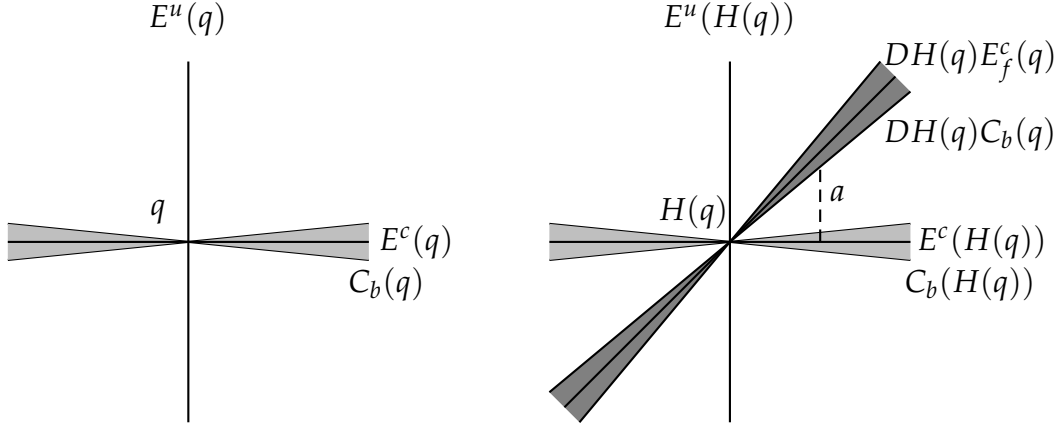

\begin{center}
\vspace{3cm}
\psset{origin={-4,0},xunit=.8cm,yunit=.8cm}
\rput(-5.5,.5){$q$}
\psline{-}(0,-3)(0,3)
\rput(-1.2,0){$E^{c}(q)$}
\rput(-5, 3.7){$E^{u}(q)$}
\psline{-}(-3,-0.3)(3,0.3)
\psline{-}(-3,0.3)(3,-0.3)
\rput(-1.7,-.7){$C_b(q)$}
\pspolygon*[linecolor=lightgray](3,0.3)(3,-0.3)
\pspolygon*[linecolor=lightgray](-3,0.3)(-3,-0.3)
\psline{-}(-3,0)(3,0)
\psset{origin={4,0},xunit=.8cm,yunit=.8cm}
\psline{-}(0,-3)(0,3)
\rput(3.4,.5){$H(q)$}
\rput(8.2,0){$E^{c}(H(q))$}
\rput(4, 3.7){$E^{u}(H(q))$}
\psline{-}(-3,-0.3)(3,0.3)
\psline{-}(-3,0.3)(3,-0.3)
\rput(8,-.7){$C_b(H(q))$}
\pspolygon*[linecolor=lightgray](3,0.3)(3,-0.3)
\pspolygon*[linecolor=lightgray](-3,0.3)(-3,-0.3)
\psline{-}(-3,0)(3,0)
\pspolygon*[linecolor=gray](3,2.5)(2.5,3)
\pspolygon*[linecolor=gray](-3,-2.5)(-2.5,-3)
\psline{-}(-3,-2.5)(3,2.5)
\psline{-}(-2.5,-3)(2.5,3)
\psline{-}(-2.75,-2.75)(2.75,2.75)
\rput(8,3){$DH(q)E^c_f(q)$}
\rput(8,1.8){$DH(q)C_b(q)$}
\psline[linestyle=dashed]{-}(1.5,0)(1.5,1.2)
\rput(5.8,.8){$a$}
\vspace{3cm}
\caption{Action of $DH(q)$ over the cone $C_b(q)$.}\label{fig:conos1}
\end{center}
\end{figure}

We can reformulate \eqref{eq:slopebounds4} in terms of
cones. For $b\geq 0$, we consider the cones inside
$E^{cu}_f(q)$ around the subspace $E^c_f(q)$ given for $q\in
M$ by
\begin{align*}
  C_b(q)&:=\{v=v^c+v^u\in  E^{cu}_f(q)\::\:
 b\ge s_c(v)\}\cup\{0\}.
\end{align*}
Then, for every $q\in B(q_0,\eta)$, $DH(q)$ carries the cone
$C_b(q)$ in a cone around $DH(q)E^c_f(q)$ whose slope is at
least $a>b>0$. In particular $DH(q)C_b(q)\cap
C_b(H(q))=\{0\}$; see Figure~\ref{fig:conos1}.

In addition, in the case $q\in B(q_0,\eta)\subseteq V$,
$g(q),\dots, g^{n}(q)\notin V$ and $g^{n+1}(q)\in V$, $n\geq
1$, we have
$$s_c(Dg^{n+1}(q)v)\geq\lambda_3^{n+1}s_c(DH(q)v)>\lambda_3^{n+1}a>s_c(v) .$$
This implies that, when $q\in B(q_0,\eta)$ returns to $V$
after $n+1\geq 2$ iterates, the cone $Dg^{n+1}C_b(q)$ is far
away from the cone $C_b(g^{n+1}(q))$. In other words, for
any $q\in B(q_0,\eta)$ the cone $C_b(g^{R(q)}(q)$ cannot be
backward $Dg^{R(q)}$-invariant, where $R(q)$ is the first
return map of $q$ to $V$ under the action of $g$. Hence, for
every $q\in B(q_0,\eta)$ returning to $V$ in a future
iterate of $g$, the central direction $E^c_g(g^{R(q)}(q))$
of $g$ cannot be contained in $C_b(g^{R(q)}(q))$.

Denote by $\tilde{V}$ the set of points $\tilde{q}\in V$
such that there are $R\geq 2$ and $q\in B(q_0,\eta)$ such
that $g^R(q)=\tilde{q}$ and for $1\leq n<R$, $g^n(q)\notin
V$. Thus, for every $\tilde{q}\in V$,
$E^c_g(g^{R(q)}(q))=E^c_g(\tilde q)$ can be written as the
graph of a nonzero linear map $G_{\tilde q}:E^c_f(\tilde
q)\to E^u_f(\tilde q)$.

As explained in Section~\ref{sec:minimal-unstable-fol-1},
our assumption of minimality of the unstable foliation
implies that $\mu$ almost every point visits any given open
subset, the set $V$ say, with positive asymptotic frequency,
for each $u$-Gibbs state $\mu$ of $g$. In particular,
$$\mu(\tilde V)=\mu\big(B(q_0,\eta)\big)>0$$
concluding the proof of the lemma.
\endproof

\section{Comparing the action of the derivatives}
\label{sec:compar-action-deriva}

Here we complete the proof of the main
Theorem~\ref{mthm:removing}, providing the details of
\eqref{eq:qoutient-gt-1} from the overview in
Section~\ref{sec:overvi-arguments}. 

We compare the norm of the actions of the derivatives of $g$
and $f$ on the new central subbundle and the old central
subbundle respectively.  We use the adapted norms,
introduced in Section~\ref{sec:adapted-norms-partia}, to our
advantage in the calculation that follow, together with the
smoothness of the strong-unstable foliation.

The main idea comes from a simple fact of linear algebra: 

\begin{lemma}\label{le:al}
  Consider $A\::\:E\oplus F\to E\oplus F$ a linear
  transformation where the subspaces $E$ and $F$ are
  invariant under the action of $A$ . Assume that there
  exists $\lambda>1$ such that
\begin{itemize}
\item[(a)] ${\rm m}(A|F)>\lambda$, that means $A$ is
  uniformly expanding on $F$; and,
\item[(b)] $\|A|E\|<\lambda$, or equivalently, the splitting
  $E\oplus F$ is dominated.
\end{itemize}
Let $L\::\:E\to F$ be a linear map, $L\not\equiv 0$, put
$G={\rm graph\:}L$ and assume that the norm above is given
by an inner product on $E\oplus F$ such that $E$ is
orthogonal to $F$.  Then there exists $\xi>0$ such that
\begin{equation}
  \|A|G\|\geq(1+\xi)\|A|E\|.
\end{equation}
\end{lemma}

\proof Fix $u\in G$, $u\ne 0$. Then, there is $v\in
E\setminus\{0\}$ such that $u=v+Lv$. Since the decomposition
$E\oplus F$ is $A$-invariant, then $v, Av\in E$ and $Lv,
ALv\in F$. Moreover,
$$
\frac{\|Au\|^2}{\|u\|^2}
=
\frac{\|Av+ALv\|^2}{\|v+Lv\|^2}
=
\frac{\|Av\|^2+\|ALv\|^2}{\|v\|^2+\|Lv\|^2}
=
\frac{\|Av\|^2}{\|v\|^2}
\left[\frac{1+\frac{\|ALv\|^2}{\|Av\|^2}}{1+\frac{\|Lv\|^2}{\|v\|^2}}\right].
$$
If $u=v+w\in E\oplus F$ with $u\in E, v\in F$ and
$s_E(u)=\frac{\|w\|}{\|v\|}$ is the slope of $u$, then it
follows directly from the assumptions over $A$ that
\begin{equation}\label{eq:leal1}s_E(Au)>s_E(u).\end{equation}
Simple algebraic manipulation over \eqref{eq:leal1} give us
\begin{equation}\label{eq:leal2}
\left[\frac{1+\frac{\|ALv\|^2}{\|Av\|^2}}{1+\frac{\|Lv\|^2}{\|v\|^2}}\right]
\geq  1+\xi
\end{equation}
and the statement of the lemma follows.
\endproof

\begin{lemma}\label{le:ML}
  There exists a measurable function $\xi\::\: M\to\RR$ such
  that for every $q\in M$,
$$\|Dg(q)|E^c_g(q)\|\geq 1+\xi(q).$$
Moreover, the set of points $q\in M$ such that $\xi(q)>0$
has positive $\mu$ measure for every $u$-Gibbs state $\mu$
of $g$.
\end{lemma}

\proof We subdivide the argument in a number of cases for
clarity.
\begin{description}
\item[CASE A] $E^c_g(q)$ is the graph of the non-zero linear
  map $G_q\::\:E^c_f(q)\to E^u_f( q)$.

In this case 
\begin{align*}
  E^c_g(q)=\{u+G_q(u): u\in E^c_f(y)\} =(I+G_q)(E^c_f(q)).
\end{align*}
\begin{description}
\item[CASE A.1]  for $q\in M\setminus V$, we can
apply the previous Lemma~\ref{le:al} to $A=Dg(q)=Df(q)$, $E=E^c_f(q)$,
$F=E^u_g(q)=E^u_f(q)$, $G=E^c_g(q)$ and $L=G_q$ obtaining
\begin{align}\label{eq:quotient2}
\|Df(q)|E^c_g(q)\|  \ge (1+\xi(q))\|Df(q)|E^c_f(q)\|=1+\xi(q)
\end{align}
for some $\xi(q)>0$, which is clearly true by our
assumptions~\eqref{PH2} and~\eqref{PH3}.
\item[CASE A.2] for $q\in V$, we apply Lemma~\ref{le:al} to
$A=Dg(q)=Df(H(q))$, $E=DH(q)E^c_f(q)$, $F=DH(q)E^u_f(q)$,
$G=DH(q)E^c_g(q)$ and $L=DH(q)\circ G_q$ obtaining
\begin{eqnarray}
  \|Dg(q)|E^c_g(q)\|&=&\|Df(H(q))|DH(q)E^c_g(q)\| \nonumber\\
  &\ge& (1+\xi(q))\|Df(H(q))|E^c_f(H(q))\|=1+\xi(q)  \label{eq:neededineq2}
\end{eqnarray}
and this is clearly true again from assumptions \eqref{PH2}
and~\eqref{PH3}.
\end{description}

\item[CASE B] $E^c_g(q)=E^c_f(q)$ (that is, we assume that $G_q\equiv0$).
  \begin{description}
  \item[CASE B.1] for $q\in M\setminus V$ we have 
    \begin{align}\label{eq:quotient3}
      \|Dg(q)|E^c_g(q)\| =\|Df(q)|E^c_f(q)\| =1.
    \end{align}
  \item[CASE B.2] otherwise, for  $q\in V$ we have
$$\|Dg(q)|E^c_g(q)\|=\|Df(H(q))|DH(q)E^c_g(q)\|
=\|Df(H(q))|DH(q)E^c_f(q)\|.$$
Now we have the following two alternatives.
\begin{description}
\item[CASE B.2a] If $DH(q)E^c_f(q)=E^c_f(q)$, then $\|Dg(q)|E^c_g(H(q))\|=1$.
\item[CASE B.2b] Otherwise, $DH(q)E^c_f(q)$ is the graph of
  a non-zero linear maps $G_{H(q)}\::\: E^c_f(H(q))\to
  E^u_q(H(q))$ so we apply Lemma~\ref{le:al} to
  $A=Dg(q)=Df(H(q))$, $E=E^c_f(H(q))$, $F=E^u_f(H(q))$,
  $G=DH(q)E^c_f(q)$ and $L= G_{H(q)}$ obtaining
\begin{eqnarray}\|Dg(q)|E^c_g(q)\|&=&\|Df(H(q))|DH(q)E^c_g(q)\| \nonumber\\
&\ge& (1+\xi(q))\|Df(H(q))|E^c_f(H(q))\|=1+\xi(q)  \label{eq:neededineq3}
\end{eqnarray}
\end{description}
\end{description}
\end{description}

Finally, putting together \eqref{eq:quotient2},
\eqref{eq:neededineq2}, \eqref{eq:quotient3} and
\eqref{eq:neededineq3} we obtain a measurable function
$\xi\::\: M\to \RR$ such that for every $q\in M$,
$\xi(q)\geq 0$ and
$$\|Dg(q)|E^c_g(q)\|\geq 1+\xi(q).$$

It follows from Lemma~\ref{le:Ecg-away-Ecf} that the set of
points $q\in M$ such that $\xi(q)>0$ has positive $\mu$
measure for every $u$-Gibbs state $\mu$ of $g$.
\endproof
 

As explained in the outline of the proof, in
Section~\ref{sec:overvi-arguments}, the main theorem and
corollaries follow easily from this estimates and known
results on mostly expanding partially hyperbolic
diffeomorphisms.
This completes the proof of the main theorem.


\def\cprime{$'$}


\end{document}